\title{$A_2$-Planar Algebras II: \\
        Planar Modules}
\author{
        David E. Evans and Mathew Pugh \\ \\
        School of Mathematics, \\
        Cardiff University, \\
        Senghennydd Road, \\
        Cardiff, CF24 4AG, \\
        Wales, U.K.
}
\date{\today}

\documentclass[12pt]{article}
\usepackage{amssymb}
\usepackage{graphicx}
\usepackage[all]{xy}

\newtheorem{Def}{Definition}[section]
\newtheorem{Prop}[Def]{Proposition}
\newtheorem{Lemma}[Def]{Lemma}
\newtheorem{Cor}[Def]{Corollary}
\newtheorem{Thm}[Def]{Theorem}

\begin{document}
\maketitle

\begin{abstract}
Generalizing Jones's notion of a planar algebra, we have previously introduced an $A_2$-planar algebra capturing the structure contained in the double complex pertaining to the subfactor for a finite $SU(3)$ $\mathcal{ADE}$ graph with a flat cell system. We now introduce the notion of modules over an $A_2$-planar algebra, and describe certain irreducible Hilbert $A_2$-$TL$-modules. We construct an $A_2$-graph planar algebra associated to each pair $(\mathcal{G},W)$ given by an $SU(3)$ $\mathcal{ADE}$ graph $\mathcal{G}$ and a cell system $W$ on $\mathcal{G}$. A partial modular decomposition of these $A_2$-graph planar algebras is achieved.
\end{abstract}

\section{Introduction}

We introduced in \cite{evans/pugh:2009iii} the notion of an $A_2$-planar algebra. This was useful to understand the double complexes of finite dimensional algebras which arise in the context of $SU(3)$ subfactors and modular invariants. Here we begin a study of their planar modules.

These $A_2$-planar algebras are a direct generalization of the planar algebras of Jones \cite{jones:planar}.
To avoid too much confusion one could refer to these planar algebras of Jones here as $A_1$-planar algebras, which naturally contain the Temperley-Lieb algebra which encodes the representation theory of quantum $SU(2)$.
Our $A_2$-planar algebras naturally encode the representation theory of quantum $SU(3)$, or in the dual Hecke picture, the finite dimensional algebras (or $A_2$-Temperley-Lieb algebras) which appear from the representations of the deformation of the symmetric group.

A braided subfactor $N \subset M$, or equivalently a module category ${}_N \mathcal{X}_N$ over the modular tensor category ${}_N \mathcal{X}_N$, yields a modular invariant via the theory of $\alpha$-induction \cite{bockenhauer/evans/kawahigashi:1999, bockenhauer/evans:2000, evans:2003}, and a non-negative integer matrix representation, or nimrep, of the Verlinde algebra realised by ${}_N \mathcal{X}_N$. In the case of $SU(3)$, ${}_N \mathcal{X}_N$ is a finite system of endomorphisms over the type $\mathrm{III}_1$ factor $N$. The classification of $SU(3)$ modular invariants was shown to be complete by Gannon \cite{gannon:1994}. Ocneanu claimed \cite{ocneanu:2000ii, ocneanu:2002} that all $SU(3)$ modular invariants were realised by subfactors and this was shown in \cite{ocneanu:2000ii, ocneanu:2002, xu:1998, bockenhauer/evans:1999i, bockenhauer/evans:1999ii, bockenhauer/evans/kawahigashi:1999, bockenhauer/evans:2001, bockenhauer/evans:2002, evans/pugh:2009i, evans/pugh:2009ii}. The braided subfactor, or module category, and its associated modular invariant, are both classified by an $SU(3)$ $\mathcal{ADE}$ graph $\mathcal{G}$, which is the graph whose adjacency matrix is given by evaluating the nimrep at the fundamental generator of ${}_N \mathcal{X}_N$.

These graphs carry a cell system \cite{ocneanu:2000ii, ocneanu:2002}. These cells give numerical weight to Kuperberg's \cite{kuperberg:1996} diagram of trivalent vertices -- corresponding to the fact that the trivial representation is contained in the triple product of the fundamental representation of $SU(3)$ through the determinant. They yield, in a natural way, representations of an $A_2$-Temperley-Lieb or Hecke algebra.
We computed the numerical values of the Ocneanu cells in \cite{evans/pugh:2009i}.
For $SU(2)$ or bipartite graphs, the corresponding weights (associated to the diagrams of cups or caps), arise in a more straightforward fashion from a Perron-Frobenius eigenvector, giving a natural representation of the Temperley-Lieb algebra or Hecke algebra.

The bipartite theory of the $SU(2)$ setting has to some degree become a three-colourable theory in our $SU(3)$ setting. This theory is not completely three-colourable since some of the graphs are not three-colourable -- namely the graphs $\mathcal{A}^{(n)\ast}$ associated to the conjugate modular invariants, $n \geq 4$, $\mathcal{D}^{(n)}$ associated to the orbifold modular invariants, $n \neq 0 \textrm{ mod } 3$, and the exceptional graph $\mathcal{E}^{(8)\ast}$. The figures for the complete list of the $\mathcal{ADE}$ graphs are given in \cite{behrend/pearce/petkova/zuber:2000, evans/pugh:2009i}.

In Section \ref{sect:JonesPA} we review the basics of Jones' planar algebras and planar modules, and in Section \ref{sect:A2planar_algebras} we review our construction of $A_2$-planar algebras. In Section \ref{sect:A2planar_modules} we describe $A_2$-planar modules, including the notion of lowest weight. In particular we describe all irreducible $A_2$-$PTL$-modules of lowest weight zero. Then in Section \ref{Sect:A2-GPA} we construct the $A_2$-graph planar algebra $P^{\mathcal{G}}$ for an $\mathcal{ADE}$ graph, and determine a partial decomposition of $P^{\mathcal{G}}$ into irreducible $A_2$-$PTL$-modules.

Before we delve into the theory of ($A_2$-)planar algebras, we review the realisation of modular invariants by braided subfactors.
Let $A$ and $B$ be type $\mathrm{III}$ von Neumann factors. A unital $\ast$-homomorphism $\rho:A\rightarrow B$ is called a $B$-$A$ morphism.
Some $B$-$A$ morphism $\rho'$ is called equivalent to $\rho$ if $\rho'=\mathrm{Ad}(u)\circ\rho$ for some unitary $u\in B$.
The equivalence class $[\rho]$ of $\rho$ is called the $B$-$A$ sector of $\rho$. If $\rho$ and $\sigma$ are $B$-$A$ morphisms with finite statistical dimensions, then the vector space of intertwiners
$\mathrm{Hom}(\rho,\sigma)=\{ t\in B: t\rho(a)=\sigma(a)t \,, \,\, a\in A \}$
is finite-dimensional, and we denote its dimension by $\langle \rho, \sigma \rangle$.
A $B$-$A$ morphism is called irreducible if $\langle \rho,\rho \rangle=1$, i.e. if $\mathrm{Hom}(\rho,\rho) = \mathbb{C} \mathbf{1}_B$.
Then, if $\langle \rho, \tau \rangle \neq 0$ for some (possibly reducible) $B$-$A$ morphism $\tau$, $[\rho]$ is called an irreducible subsector of $[\tau]$ with multiplicity $\langle \rho, \tau \rangle$.
An irreducible $A$-$B$ morphism $\overline{\rho}$ is a conjugate morphism of the irreducible $\rho$ if and only if $[\overline{\rho}\rho]$ contains the trivial sector $[\mathrm{id}_A]$ as a subsector, and then $\langle \rho\overline{\rho}, \mathrm{id}_B \rangle = 1 = \langle \overline{\rho}\rho, \mathrm{id}_A \rangle$ automatically \cite{izumi:1991}.

The Verlinde algebra is realised in the subfactor models by systems of endomorphisms ${}_N \mathcal{X}_N$ of the hyperfinite type $\mathrm{III}_1$ factor $N$.
That is, ${}_N \mathcal{X}_N$ denotes a finite system of finite index irreducible endomorphisms of a factor $N$ such that the elements of ${}_N \mathcal{X}_N$ are not unitary equivalent, for any $\lambda \in {}_N \mathcal{X}_N$ there is a representative $\overline{\lambda} \in {}_N \mathcal{X}_N$ of the conjugate sector $[\overline{\lambda}]$, and ${}_N \mathcal{X}_N$ is closed under composition and subsequent irreducible decomposition.
In the case of WZW models associated to $SU(n)$ at level $k$, the Verlinde algebra is a non-degenerately braided system of endomorphisms ${}_N \mathcal{X}_N$, labelled by the positive energy representations of the loop group of $SU(n)_k$ on a type $\mathrm{III}_1$ factor $N$, with fusion rules $\lambda \mu = \bigoplus_{\nu} N_{\lambda \nu}^{\mu} \nu$ which exactly match those of the positive energy representations \cite{wassermann:1998}. The fusion matrices $N_{\lambda} = [N_{\rho \lambda}^{\sigma}]_{\rho,\sigma}$ are a family of commuting normal matrices which give a representation themselves of the fusion rules of the positive energy representations of the loop group of $SU(n)_k$, $N_{\lambda} N_{\mu} = \sum_{\nu} N_{\lambda \nu}^{\mu} N_{\nu}$.
This family $\{ N_{\lambda} \}$ of fusion matrices can be simultaneously diagonalised:
\begin{equation} \label{eqn:verlinde_formula}
N_{\lambda} = \sum_{\sigma} \frac{S_{\sigma, \lambda}}{S_{\sigma,1}} S_{\sigma} S_{\sigma}^{\ast},
\end{equation}
where $1$ is the trivial representation, and the eigenvalues $S_{\sigma, \lambda}/S_{\sigma,1}$ and eigenvectors $S_{\sigma} = [S_{\sigma, \mu}]_{\mu}$ are described by the statistics $S$ matrix.
Moreover, there is equality between the statistics $S$- and $T$- matrices and the Kac-Peterson modular $S$- and $T$- matrices which perform the conformal character transformations \cite{kac:1990}, thanks to \cite{frohlich/gabbiani:1990, fredenhagen/rehren/schroer:1992, wassermann:1998}.

The key structure in the conformal field theory is the modular invariant partition function $Z$. In the subfactor setting this is realised by
a braided subfactor $N \subset M$ where trivial (or permutation) invariants in the ambient factor $M$ when restricted to $N$ yield $Z$. This would mean that the dual canonical endomorphism decomposes as a finite linear combination of endomorphisms in ${}_N \mathcal{X}_N$.
Indeed if this is the case for the inclusion $N \subset M$, then the process of $\alpha$-induction allows us to analyse the modular invariant,
providing two extensions of $\lambda$ on $N$ to endomorphisms $\alpha^{\pm}_{\lambda}$ of $M$, such that the matrix $Z_{\lambda,\mu} = \langle \alpha_{\lambda}^+, \alpha_{\mu}^- \rangle$ is a modular invariant \cite{bockenhauer/evans/kawahigashi:1999, bockenhauer/evans:2000, evans:2003}.

Let ${}_N \mathcal{X}_M$, ${}_M \mathcal{X}_M$ denote a system of endomorphisms consisting of a choice of representative endomorphism of each irreducible subsector of sectors of the form $[\lambda \overline{\iota}]$, $[\iota \lambda \overline{\iota}]$ respectively, for each $\lambda \in {}_N \mathcal{X}_N$, where $\iota: N \hookrightarrow M$ is the inclusion map which we may consider as an $M$-$N$ morphism, and $\overline{\iota}$ is a representative of its conjugate $N$-$M$ sector.
The action of the system  ${}_N \mathcal{X}_N$ on the $N$-$M$ sectors ${}_N \mathcal{X}_M$ produces a nimrep (non-negative matrix integer representation of the fusion rules) $G_{\lambda} G_{\mu} = \sum_{\nu} N_{\lambda \nu}^{\mu} G_{\nu}$,
whose spectrum reproduces exactly the diagonal part of the modular invariant, i.e.
\begin{equation} \label{eqn:verlinde_formulaG}
G_{\lambda} = \sum_{\sigma} \frac{S_{\sigma,\lambda}}{S_{\sigma,1}} \psi_{\sigma} \psi_{\sigma}^{\ast},
\end{equation}
with the spectrum of $G_{\lambda} = \{ S_{\mu, \lambda}/S_{\mu,1}$ with multiplicity $Z_{\mu,\mu} \}$ \cite[Theorem 4.16]{bockenhauer/evans/kawahigashi:2000}. The labels $\mu$ of the non-zero diagonal elements are called the exponents of $Z$, counting multiplicity. A modular invariant for which there exists a nimrep whose spectrum is described by the diagonal part of the invariant is said to be nimble.

The systems ${}_N \mathcal{X}_N$, ${}_N \mathcal{X}_M$, ${}_M \mathcal{X}_M$ are (the irreducible objects of) tensor categories of endomorphisms with the Hom-spaces as their morphisms. Thus ${}_N \mathcal{X}_N$ gives a braided modular tensor category, and ${}_N \mathcal{X}_M$ a module category.
The structure of the module category ${}_N \mathcal{X}_M$ is the same as a tensor functor $F$ from ${}_N \mathcal{X}_N$ to the category $\mathrm{Fun}({}_N \mathcal{X}_M,{}_N \mathcal{X}_M)$ of additive functors from ${}_N \mathcal{X}_M$ to itself, see \cite{ostrik:2003}.

In \cite{evans/pugh:2010ii} we described a $\mathbb{C}$-linear tensor category $TL$ called the Temperley-Lieb category, which depends on a paramater $\delta = q+q^{-1}$, where $q$ is real or a root of unity. Its objects are given by direct sums of Jones-Wenzl projections \cite{wenzl:1987}, whilst its morphisms are matrices whose entries are planar tangles. In the language of Section \ref{sect:JonesPA} these planar tangles are planar $j$-tangles, $j \geq 0$, with no internal discs and with the outer disc removed. The Temperley-Lieb category is semisimple for $q$ real, i.e. $2 \geq \delta \in \mathbb{R}$. For $q$ a $m^{\mathrm{th}}$ root of unity, the quotient $TL^{(m)} := TL/I^{(m)}$ is semisimple, where $I^{(m)}$ is the unique proper tensor ideal in $TL$ generated by the morphism $\frak{f}_m = \mathrm{id}_{f_m}$, and $f_m$ is the $m^{\mathrm{th}}$ Jones-Wenzl projection.

For $q$ a $k+2^{\mathrm{th}}$ root of unity, these Jones projections satisfy the same fusion relations as the positive energy representations of the loop group of $SU(2)_{k}$, and hence each object of $TL^{(k+2)}$ can be identified naturally with an endomorphism in ${}_N \mathcal{X}_N$.
Then for a braided subfactor $N \subset M$ with classifying graph $\mathcal{G}$ one of the $ADE$ graphs, that is $N \subset M$ yields nimrep $G$ such that $G_{\rho} = \mathcal{G}$, ${}_N \mathcal{X}_M$ defines a tensor functor $F$ from $TL^{(k+2)}$ to $\mathrm{Fun}({}_N \mathcal{X}_M,{}_N \mathcal{X}_M)$. This functor recovers the bipartite graph planar algebra construction \cite{jones:2000} described here in Section \ref{Sect:planar_alg_of_bipartite}.

A parallel construction of the $A_2$-Temperley-Lieb category $A_2$-$TL$ was given in \cite{evans/pugh:2010ii}.
Then for a braided subfactor $N \subset M$ with classifying graph $\mathcal{G}$ one of the $\mathcal{ADE}$ graphs, which labels a modular invariant of $SU(3)$ at level $k$, ${}_N \mathcal{X}_M$ defines a tensor functor $F: A_2\mathrm{-}TL^{(k+3)} \longrightarrow \mathrm{Fun}({}_N \mathcal{X}_M,{}_N \mathcal{X}_M)$ which recovers the construction of the $A_2$-planar algebra for an $SU(3)$ $\mathcal{ADE}$ graph. This construction is described here in Section \ref{Sect:A2-GPA}.

It would be interesting to extend the $A_2$-planar algebra formalism to other rank 2 simple Lie algebras such as $G_2$ or $C_2 = sp(4)$, based on the web spaces of Kuperberg \cite{kuperberg:1996}. Generalized Jones-Wenzl idempotents for $C_2$ were constructed in \cite{kim:2007}. One could also try to extend the formalism to $A_n$, $n \geq 3$.
A complete list of relations for a diagrammatic formalism for $A_3$ has been conjectured \cite{kim:2003}. The generators for this formalism yield cells on the classifying graphs for $SU(4)$ modular invariants, and the relations yield consistency relations for these cells. These graphs were proposed by Ocneanu \cite{ocneanu:2002}.
In order to construct the corresponding ($A_3$-)planar algebra, one would need to compute the values of these cell systems.
Some $SU(4)$ modular invariants have been realised using braided subfactors which come from conformal and orbifold embeddings \cite{xu:1998, bockenhauer/evans:1999i, bockenhauer/evans:1999ii, ocneanu:2002}.

\section{Preliminaries on Jones' planar algebras and planar modules} \label{sect:JonesPA}

Let us briefly review the basics of Jones' planar algebras, and the notion of planar modules over these algebras \cite{jones:planar,jones:2001}. To avoid confusion we will refer to these planar algebras and modules as $A_1$-planar algebras and modules. A planar $k$-tangle consists of a disc $D$ in the plane with $2k$ vertices on its boundary, $k \geq 0$, and $n \geq 0$ internal discs $D_j$, $j=1,\ldots,n$, where the disc $D_j$ has $2k_j$ vertices on its boundary, $k_j \geq 0$. One vertex on the boundary of each disc (including the outer disc $D$) is chosen as a marked vertex, and the segment of the boundary of each disc between the marked vertex and the vertex immediately adjacent to it as we move around the boundary in an anti-clockwise direction is labelled $\ast$. Inside $D$ we have a collection of disjoint smooth curves, called strings, where any string is either a closed loop, or else has as its endpoints the vertices on the discs, and such that every vertex is the endpoint of exactly one string. Any tangle must also allow a checkerboard colouring of the regions inside $D$.

The planar operad is the collection of all diffeomorphism classes of such planar tangles, with composition of planar tangles defined.
A planar algebra $P$ is then defined to be an algebra over this operad, i.e. a family $P = (P_k^+, P_k^-; k \geq 0)$ of vector spaces with $P_k^{\pm} \subset P_{k'}^{\pm}$ for $k < k'$, and with the following property: for every $k$-tangle $T$ with $n$ internal discs $D_j$ labelled by elements $x_j \in P_{k_j}$, $j=1,\ldots,n$, there is an associated linear map $Z(T):\otimes_{j=1}^n P_{k_j} \rightarrow P_k$, which is compatible with the composition of tangles and re-ordering of internal discs.

A planar module over $P$ is a graded vector space $V = (V_k^+, V_k^-; k \geq 0)$ with an action of $P$. Given a planar $m$-tangle $T$ in $\mathcal{P}$ with distinguished ($V$ input) internal disc $D_1$ with with $2k$ vertices on its boundary, $k \geq 0$, and other ($P$ input) internal discs $D_p$, $p=2, \ldots, n$, with $2k_p$ vertices on its boundary, $k_p \geq 0$, there is a linear map $Z(T): V_k \otimes \left( \otimes_{p=2}^n P_{k_p} \right) \rightarrow V_m$, where $Z(T)$ satisfies the same compatability conditions as for $P$ itself.

\subsection{$P^{\mathcal{G}}$ as a $TL$-module for an $ADE$ Dynkin diagram $\mathcal{G}$} \label{Sect:planar_alg_of_bipartite}

Jones \cite{jones:2001} determined all Hilbert Temperley-Lieb modules $H^{k,\omega}$ of lowest weight $k > 0$, $k \in \mathbb{N}$, and $H^{\mu}$ of lowest weight 0. We review these modules. For $k,m \in \mathbb{N}$, let $ATL_{m,k}$ denote the space of all annular $(m,k)$-tangles (having $2m$ vertices on the outer disc and $2k$ vertices on the (distinguished) inner disc, where the vertices have alternating orientations) with no other internal discs. Tangles are composed by inserting one annular $(m,k)$-tangle inside the internal disc of an annular $(k,n)$-tangle. For $1 \leq k \leq m$, $k,m \in \mathbb{N}$, let $\mathcal{T}_m^k$ denote the set of annular $(m,k)$-tangles with no internal discs and $2k$ through strings. If $\widetilde{ATL}_{m,k}$ denotes the quotient of $ATL_{m,k}$ by the ideal generated by all annular $(m,k)$-tangles with no internal discs and strictly less than $2k$ through strings, then the equivalence classes of the elements of $\mathcal{T}_m^k$ form a basis for $\widetilde{ATL}_{m,k}$. The group $\mathbb{Z}_k$ acts by an internal rotation, which permutes the basis elements. The action of $ATL$ on $\widetilde{ATL}_{m,k}$ is given as follows. Let $T$ be an annular $(p,m)$-tangle in $ATL_{p,m}$ and $R \in \mathcal{T}_m^k$. Define $T(R)$ to be $\delta^r \widehat{TR}$ if the $(p,k)$-tangle $TR$ has $2k$ through strings and 0 otherwise, where $TR$ contains $r$ contractible circles and $\widehat{TR}$ is the tangle $TR$ with all the contractible circles removed. Since the action of $ATL$ commutes with the action of $\mathbb{Z}_k$, $\widetilde{ATL}_{m,k}$ splits as a $TL$-module into a direct sum, over the $k^{\textrm{th}}$ roots of unity $\omega$, of $TL$-modules $V^{k,\omega}_m$ which are the eigenspaces for the action of $\mathbb{Z}_k$ with eigenvalue $\omega$. For each $k$ one can choose a faithful trace $\mathrm{tr}$ on the abelian $C^{\ast}$-algebra $\widetilde{ATL}_{k,k}$, which extends to $ATL_{k,k}$ by composition with the quotient map. The inner-product on $\widetilde{ATL}_{m,k}$ is then defined to be $\langle S, T \rangle = \mathrm{tr}(T^{\ast}S)$ for $S,T \in \widetilde{ATL}_{m,k}$.

We now turn to the zero-weight case ($k=0$). The algebras $ATL_{\pm}$, which have the regions adjacent to both inner and outer boundaries shaded $\pm$, are generated by elements $\sigma_{\pm} \sigma_{\mp}$, where $\sigma_{\pm}$ is the $(\pm,\mp)$-tangle which is just a single non-contractible circle, with the region which meets the outer boundary shaded $\pm$ and the region which meets the inner boundary shaded $\mp$. Then the dimensions on $V_+$ and $V_-$ must be 1 or 0 for any $TL$-module $V$. In $V$, the maps $\sigma_{\pm} \sigma_{\mp}$ must contribute a scalar factor $\mu^2$, where $0 \leq \mu \leq \delta$. If $\mu = \delta$, $V^{\delta}$ is simply the ordinary Temperley-Lieb algebra.
When $0 < \mu < \delta$, $V^{\mu}$ is the $TL$-module such that $V^{\mu}_m$ $m \geq 0$, has as basis the set of $(m,+)$-tangles with no internal discs and at most one non-contractible circle. The action of $ATL$ on $V^{\mu}$, $0 \leq \mu \leq \delta$, is given as follows. Let $T$ be an annular $(p,m)$-tangle in $ATL_{p,m}$ and $R$ be a basis element of $V^{\mu}$. Define $T(R)$ to be $\delta^r \mu^{2d} \widehat{TR}$, where $TR$ contains $r$ contractible circles and $2d+i$ non-contractible circles, where $i \in \{ 0,1 \}$, and $\widehat{TR}$ is the tangle $TR$ with all the contractible circles removed and $2d$ of the non-contractible circles removed. The inner product on $V^{\mu}$ is defined by $\langle S,T \rangle = \delta^r \mu^{2d}$, where $T^{\ast}S$ contains $r$ contractible circles and $2d$ non-contractible circles.
When $\mu = 0$, we have $TL$-modules $V^{0,+}$ and $V^{0,-}$, where $V^{0,\pm}_m$ has as basis the set of $(m,\pm)$-tangles with no internal discs and no contractible circles. The action of $ATL$ on $V^{0,\pm}$ is given as follows. Let $T$ be an annular $(p,m)$-tangle in $ATL_{p,m}$ and $R$ be a basis element of $V^{0,\pm}$. Define $T(R)$ to be $\delta^r \widehat{TR}$, where $TR$ contains $r$ contractible circles. Now $\widehat{TR}$ is zero if $TR$ contains any non-contractible circles, and is the tangle $TR$ with all the contractible circles removed otherwise. The inner product on $V^{0,\pm}$ is defined by $\langle S,T \rangle = 0$ if $T^{\ast}S$ contains any non-contractible circles, and $\langle S,T \rangle = \delta^r$ otherwise, where $r$ is the number of contractible circles in $T^{\ast}S$.

In the generic case, $\delta > 2$, it was shown \cite{jones:2001} that the inner-product is always positive definite, so that $H = V$ is a Hilbert $TL$-module, for the irreducible lowest weight $TL$-module $V$. In the non-generic case, if the inner product is positive semi-definite, $H$ is defined to be the quotient of $V$ by the vectors of zero-length with respect to the inner product.

Let $\mathcal{G}$ be a bipartite graph. Then the vertex set of $\mathcal{G}$ is given by $\mathfrak{V} = \mathfrak{V}_+ \cup \mathfrak{V}_-$, where each edge connects a vertex in $\mathfrak{V}_+$ to one in $\mathfrak{V}_-$. We call the vertices in $\mathfrak{V}_+$, $\mathfrak{V}_-$ the even, odd respectively vertices of $\mathcal{G}$. We will use the convention that the distinguished vertex $\ast$ of $\mathcal{G}$, which has the highest Perron-Frobenius weight, is an even vertex. The adjacency matrix of $\mathcal{G}$ can thus be written in the form $\left( \begin{array}{cc} 0 & \Lambda_{\mathcal{G}} \\ \Lambda_{\mathcal{G}}^T & 0 \end{array} \right)$. We let $r_{\pm} = |\mathfrak{V}_{\pm}|$.
The bipartite graph planar algebra $P^{\mathcal{G}}$ of the graph $\mathcal{G}$ was constructed in \cite{jones:2000}, which is the path algebra on $\mathcal{G}$ where paths may start at any of the even vertices of $\mathcal{G}$, and where the $m^{\textrm{th}}$ graded part $P^{\mathcal{G}}_m$ is given by all pairs of paths of length $m$ on $\mathcal{G}$ which start at the same even vertex and have the same end vertex. Let $\mu_j$, $j=1,\ldots,r_+$, denote the eigenvalues of $\Lambda_{\mathcal{G}} \Lambda_{\mathcal{G}}^T$. Then the following result is given in \cite[Prop. 13]{reznikoff:2005}: The irreducible weight-zero submodules of $P^{\mathcal{G}}$ are $H^{\mu_j}$, $j=1,\ldots,r_-$, and $r_+ - r_-$ copies of $H^0$.

Reznikoff \cite{reznikoff:2005} computed the decomposition of $P^{\mathcal{G}}$ as a $TL$-module into irreducible $TL$-modules for the $ADE$ Dynkin diagrams. For the graphs $A_m$, $m \geq 3$,
\begin{equation}
P^{A_m} = \bigoplus_{j=1}^s H^{\mu_j},
\end{equation}
where $s = \lfloor (m+1)/2 \rfloor$ is the number of even vertices of $A_m$ and $\mu_j = 2 \cos (j \pi/(m+1))$, $j=1,\ldots,s$. For $D_m$, $m \geq 3$,
\begin{equation}
P^{D_m} = \bigoplus_{j=1}^t H^{\mu_j} \oplus (s-t) H^{0,\pm} \oplus \bigoplus_{j=1}^{s-2} H^{2j,-1},
\end{equation}
where $s = \lfloor (m+2)/2 \rfloor$, $t = \lfloor (m-1)/2 \rfloor$ are the number of even, odd vertices respectively of $D_m$, and $\mu_j = 2 \cos ((2j-1) \pi/(2m-2))$, $j=1,\ldots,t$. For the exceptional graphs the results are
\begin{eqnarray}
P^{E_6} & = & H^{\mu_1} \oplus H^{\mu_4} \oplus H^{\mu_5} \oplus H^{2,-1} \oplus H^{3,\omega} \oplus H^{3,\omega^{-1}}, \\
P^{E_7} & = & H^{0,\pm} \oplus H^{\mu_1} \oplus H^{\mu_5} \oplus H^{\mu_7} \oplus H^{2,-1} \oplus H^{3,\omega} \oplus H^{3,\omega^{-1}} \oplus H^{4,-1} \oplus H^{8,-1}, \\
P^{E_8} & = & H^{\mu_1} \oplus H^{\mu_7} \oplus H^{\mu_{11}} \oplus H^{\mu_{13}} \oplus H^{2,-1} \oplus H^{3,\omega} \oplus H^{3,\omega^{-1}} \oplus H^{4,-1} \nonumber \\
& & \oplus H^{5,\zeta} \oplus H^{5,\zeta^{-1}} \oplus H^{5,\zeta^2} \oplus H^{5,\zeta^{-2}},
\end{eqnarray}
where $\omega = e^{2 \pi i/3}$, $\zeta = e^{2 \pi i/5}$, and $\mu_j = 2 \cos (\pi j/h)$ where $h$ is the Coxeter number.

\section{$A_2$-Planar Algebras} \label{sect:A2planar_algebras}

\begin{figure}[bt]
\begin{center}
  \includegraphics[width=40mm]{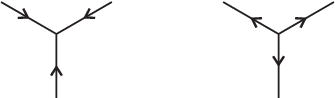}\\
 \caption{Trivalent vertices}\label{fig:A2-webs}
\end{center}
\end{figure}

We will now review the basics of $A_2$-planar algebras from \cite{evans/pugh:2009iii}.
Let $\sigma = \sigma_1 \cdots \sigma_m$ be a sign string, $\sigma_j \in \{ \pm \}$, such that the difference between the number of `$+$' and `$-$' is 0 mod 3. An $A_2$-planar $\sigma$-tangle will be the unit disc $D=D_0$ in $\mathbb{C}$ together with a finite (possibly empty) set of disjoint sub-discs $D_1, D_2, \ldots, D_n$ in the interior of $D$. Each disc $D_k$, $k \geq 0$, will have $m_k \geq 0$ vertices on its boundary $\partial D_k$, whose orientations are determined by sign strings $\sigma^{(k)} = \sigma^{(k)}_1 \cdots \sigma^{(k)}_{m_k}$ where `$+$' denotes a sink and `$-$' a source. The disc $D_k$ will be said to have pattern $\sigma^{(k)}$.
Inside $D$ we have an $A_2$-tangle where the endpoint of any string is either a trivalent vertex (see Figure \ref{fig:A2-webs}) or one of the vertices on the boundary of a disc $D_k$, $k=0, \ldots, n$, or else the string forms a closed loop. Each vertex on the boundaries of the $D_k$ is the endpoint of exactly one string, which meets $\partial D_k$ transversally.

The regions inside $D$ have as boundaries segments of the $\partial D_k$ or the strings. These regions are labelled $\overline{0}$, $\overline{1}$ or $\overline{2}$, called the colouring, such that if we pass from a region $R$ of colour $\overline{a}$ to an adjacent region $R'$ by passing to the right over a vertical string with downwards orientation, then $R'$ has colour $\overline{a+1}$ (mod 3). We mark the segment of each $\partial D_k$ between the last and first vertices with $\ast_{b_k}$, $b_k \in \{0,1,2\}$, so that the region inside $D$ which meets $\partial D_k$ at this segment is of colour $\overline{b_k}$, and the choice of these $\ast_{b_k}$ must give a consistent colouring of the regions.
For each $\sigma$ we have three types of tangle, depending on the colour $\overline{b}$ of the marked segment, or of the marked region near $\partial D$ for $\sigma = \varnothing$.

An $A_2$-planar $\sigma$-tangle $T$ with an internal disc $D_l$ with pattern $\sigma_l = \sigma'$ can be composed with an $A_2$-planar $\sigma'$-tangle $S$ with external disc $D'$ and $\ast_{D'}=\ast_{D_l}$, giving a new $\sigma$-tangle $T \circ_l S$, by inserting the $A_2$-tangle $S$ inside the inner disc $D_l$ of $T$ such that the vertices on the outer disc of $S$ coincide with those on the disc $D_l$ and the regions marked by $\ast$ also coincide. The boundary of the disc $D_l$ is removed, and the strings smoothed if necessary.
Let $\widetilde{\mathcal{P}}$ be the collection of all diffeomorphism classes of such $A_2$-planar tangles, with composition defined as above.
The $A_2$-planar operad $\mathcal{P}$ is the quotient of $\widetilde{\mathcal{P}}$ by the Kuperberg relations K1-K3 below, which are relations on a local part of the diagram: \\
\vspace{-10mm}
\begin{center}
\begin{minipage}[b]{11.5cm}
 \begin{minipage}[t]{3cm}
  \parbox[t]{2cm}{\begin{eqnarray*}\textrm{K1:}\end{eqnarray*}}
 \end{minipage}
 \begin{minipage}[t]{5.5cm}
  \begin{center}
  \mbox{} \\
 \includegraphics[width=20mm]{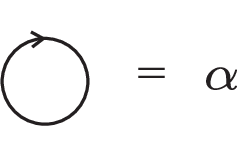}
  \end{center}
 \end{minipage}
 \begin{minipage}[t]{2cm}
  \mbox{} \\
  \parbox[t]{1cm}{}
 \end{minipage}
\end{minipage}
\begin{minipage}[b]{11.5cm}
 \begin{minipage}[t]{3cm}
  \parbox[t]{2cm}{\begin{eqnarray*}\textrm{K2:}\end{eqnarray*}}
 \end{minipage}
 \begin{minipage}[t]{5.5cm}
  \begin{center}
  \mbox{} \\
 \includegraphics[width=23mm]{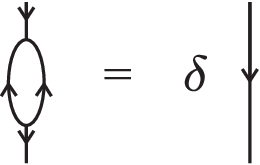}
  \end{center}
 \end{minipage}
 \begin{minipage}[t]{2cm}
  \mbox{} \\
  \parbox[t]{1cm}{}
 \end{minipage}
\end{minipage}
\begin{minipage}[b]{11.5cm}
 \begin{minipage}[t]{3cm}
  \parbox[t]{2cm}{\begin{eqnarray*}\textrm{K3:}\end{eqnarray*}}
 \end{minipage}
 \begin{minipage}[t]{5.5cm}
  \begin{center}
  \mbox{} \\
 \includegraphics[width=55mm]{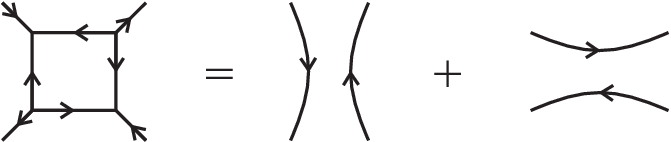}
  \end{center}
 \end{minipage}
 \begin{minipage}[t]{2cm}
  \mbox{} \\
  \parbox[t]{1cm}{}
 \end{minipage}
\end{minipage}
\end{center}
where $\delta = [2]_q$, $\alpha = [3]_q = \delta^2 - 1$, and the quantum number $[m]_q$ is defined by $[m]_q = (q^m - q^{-m})/(q-q^{-1})$, for some variable $q \in \mathbb{C}$. We will call the local picture $\,$ \includegraphics[width=10mm]{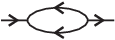} $\,$ a digon, and $\,$ \includegraphics[width=8mm]{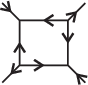} $\,$ an embedded square.

An $A_2$-planar algebra is then defined to be an algebra over this operad, i.e. a family
$P = \left( P_{\sigma}^{\overline{a}}| \, \textrm{sign strings } \sigma, a \in \{0,1,2\} \right)$
of vector spaces
with the following property: for every $\sigma$-tangle $T \in \mathcal{P}_{\sigma}$ with outer disc marked by $\ast_b$, and with $n$ internal discs $D_j$ pattern $\sigma_k$, outer disc marked by $\ast_{b_k}$ and labelled by elements $x_j \in P_{k_j}$, $j=1,\ldots,n$, there is associated a linear map $Z(T): \otimes_{k=1}^n P_{\sigma_k}^{\overline{b_k}} \longrightarrow P_{\sigma}^{\overline{b}}$ which is compatible with the composition of tangles in the following way. If $S$ is a $\sigma_k$-tangle with internal discs $D_{n+1}, \ldots, D_{n+m}$, where $D_k$ has pattern $\sigma_k$, then the composite tangle $T \circ_l S$ is a $\sigma$-tangle with $n+m-1$ internal discs $D_k$, $k = 1,2, \ldots l-1, l+1, l+2, \ldots, n+m$. From the definition of an operad, associativity means that the following diagram commutes:
\begin{equation} \label{eqn:compatability_condition_for_Z(T)}
\xymatrix{
{\left( \bigotimes_{\stackrel{k=1}{\scriptscriptstyle{k \neq l}}}^n P_{\sigma_k}^{\overline{b_k}} \right) \otimes \left( \bigotimes_{k=n+1}^{n+m} P_{\sigma_k}^{\overline{b_k}} \right)} \ar[d]_{\mathrm{id} \otimes Z(S)} \ar[dr]^(.6){Z(T \circ_l S)} \\
{\bigotimes_{k=1}^n P_{\sigma_k}}^{\overline{b_k}} \ar[r]_{Z(T)} & P_{\sigma}^{\overline{b}} }
\end{equation}
so that $Z(T \circ_l S) = Z(T')$, where $T'$ is the tangle $T$ with $Z(S)$ used as the label for disc $D_l$. We also require $Z(T)$ to be independent of the ordering of the internal discs, that is, independent of the order in which we insert the labels into the discs.
When $\sigma = \varnothing$, we will often write $P_{\varnothing}^{\overline{a}}$ as $P_{\overline{a}}$, and we adopt the convention that the empty tensor product is the complex numbers $\mathbb{C}$.

Let $\sigma^{\ast}$ be the sign string obtained by reversing the string $\sigma$ and flipping all its signs.
When each $P_{\sigma\sigma^{\ast}}$ is a $\ast$-algebra, the adjoint $T^{\ast} \in \mathcal{P}_{\sigma^{\ast}}$ of a tangle $T \in \mathcal{P}_{\sigma}$ is defined by reflecting the whole tangle about the horizontal line that passes through its centre and reversing all orientations. The labels $x_k \in P_{\sigma_k}$ of $T$ are replaced by labels $x_k^{\ast}$ in $T^{\ast}$. For any linear combination of tangles in $\mathcal{P}_{\sigma}$ the involution is the conjugate linear extension.
Then $P$ is an $A_2$-planar $\ast$-algebra if each $P_{\sigma\sigma^{\ast}}$ is a $\ast$-algebra, and for a $\sigma$-tangle $T$ with internal discs $D_k$ with patterns $\sigma_k$, labelled by $x_k \in P_{\sigma_k}$, $k=1, \ldots, n$, we have
$Z(T)^{\ast} = Z(T^{\ast})$,
where the labels of the discs in $T^{\ast}$ are $x_k^{\ast}$, and where the definition of $Z(T)^{\ast}$ is extended to linear combinations of $\sigma$-tangles by conjugate linearity.
Note a typographical error in the definition of an $A_2$-planar $\ast$-algebra in \cite[Section 4.4]{evans/pugh:2009iii}, with $P_{\sigma\sigma^{\ast}}$ above incorrectly given as $P_{\sigma}$.

\section{$A_2$-Planar Modules and $A_2$-ATL} \label{sect:A2planar_modules}

We extend Jones's notion of planar algebra modules and the annular Temperley algebra to our $A_2$-planar algebras, parallel to \cite[Section 2]{jones:2001}.

An \textbf{$A_2$-annular tangle} $T$ will be a tangle in $\mathcal{P}$ with the choice of a distinguished internal disc, which we will call the inner disc. In particular, $T$ will be called an \textbf{$A_2$-annular $(\sigma,\sigma')$-tangle} if it is an $A_2$-annular tangle with pattern $\sigma$ on its outer disc and pattern $\sigma'$ on its inner disc. If $\sigma = \varnothing$ or $\sigma' = \varnothing$, we replace the $\varnothing$ with $\overline{a}$, $a \in \{ 0,1,2 \}$, corresponding to the colour of the region which meets the outer or inner disc respectively. When $\sigma = \sigma'$ we will call $T$ an $A_2$-annular $\sigma$-tangle.
This notion of annular tangle is different to those defined in \cite{evans/pugh:2009iii}- here more than one internal disc is allowed, but one of those is chosen to be the distinguished disc (the inner boundary of the annulus).
If $P$ is an $A_2$-planar algebra, a \textbf{module} over $P$, or \textbf{$P$-module}, will be a graded vector space $V=(V_{\sigma,\overline{a}}| \, \textrm{sign strings } \sigma, a \in \{0,1,2\})$ with an action of $P$. Given an $A_2$-annular $(\sigma,\sigma')$-tangle $T$ in $\mathcal{P}$ with outer disc marked by $\ast_a$, with a distinguished ($V$ input) internal disc $D_1$ marked by $\ast_{a'}$ and with pattern $\sigma'$, and with other ($P$ input) internal discs $D_p$, $p=2, \ldots, n$, marked by $\ast_{a_p}$, with patterns $\sigma^{(p)}$, there is a linear map $Z(T): V_{\sigma',\overline{a'}} \otimes \left( \otimes_{p=2}^n P^{\overline{a_p}}_{\sigma^{(p)}} \right) \rightarrow V_{\sigma,\overline{a}}$. The map $Z(T)$ satisfies the same compatability condition (\ref{eqn:compatability_condition_for_Z(T)}) for the composition of tangles as $P$ itself.
We write $V_{\sigma}$ for the $\sigma$-graded part of $V$: $V_{\sigma} = \bigoplus_{a=0}^{2} V_{\sigma,\overline{a}}$.
An $A_2$-planar algebra is always a module over itself, called the \textbf{trivial module}.
Any relation (i.e. linear combination of labelled $A_2$-planar tangles) that holds in $P$ will hold in $V$, e.g. K1-K3 hold in $V$ where $\alpha, \delta$ have the same values as in $P$.

A module over an $A_2$-planar algebra $P$ can be understood as a module over the $A_2$-annular algebra $A_2$-$AP$, which is defined as follows. We define the associated annular category $A_2$-$AnnP$ to have objects $(\sigma,\overline{a})$, where $\sigma$ is any sign string and $a \in \{ 0,1,2 \}$, and whose morphisms are $A_2$-annular labelled tangles with labelling set all of $P$. A morphism from $(\sigma,\overline{a})$ to $(\sigma',\overline{a'})$ will be an annular $(\sigma,\sigma')$-tangle with outer disc marked by $\ast_{a}$ and inner disc by $\ast_{a'}$. Let $A_2$-$FAP$ be the linearization of $A_2$-$AnnP$- it has the same objects, but the set of morphisms from object $(\sigma,\overline{a})$ to object $(\sigma',\overline{a'})$ is the vector space having as basis the morphisms in $A_2$-$AnnP$ from $(\sigma,\overline{a})$ to $(\sigma',\overline{a'})$. Composition of morphisms is $A_2$-$FAP$ is by linear extension of composition in $A_2$-$AnnP$. The \textbf{$A_2$-annular algebra} $A_2$-$AP = \{A_2$-$AP((\sigma,\overline{a}),(\sigma',\overline{a'})) \}$ is the quotient of $A_2$-$FAP$ by relations K1-K3.
We define $A_2$-$AP_{\sigma}$ to be the algebra $\{ A_2\textrm{-}AP((\sigma,\overline{a}),(\sigma,\overline{a'}))| \, a,a' \in \{ 0,1,2 \} \}$, and define $A_2$-$AP_{\sigma}^{\overline{a}}$ to be the algebra $A_2\textrm{-}AP((\sigma,\overline{a}),(\sigma,\overline{a}))$.

\begin{figure}[tb]
\begin{center}
  \includegraphics[width=40mm]{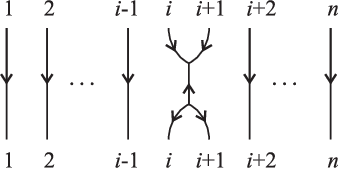}\\
 \caption{The $m$-tangle $W_i$, $i=1,\ldots, m-1$.}\label{fig:W_i}
\end{center}
\end{figure}

For for fixed $\delta \in \mathbb{C}$, we apply this procedure to the $A_2$-planar algebra $A_2$-$PTL$ defined in \cite{evans/pugh:2009iii}, which contains the algebra $\mathcal{V}^{A_2} = \mathrm{alg}(1, W_i, i \geq 1)$ as a subalgebra, where $W_i$ are the tangles illustrated in Figure \ref{fig:W_i}.
It was shown in \cite{evans/pugh:2009iii} that for real $\delta >2$, $\mathcal{V}^{A_2}$ is isomorphic to the $A_2$-Temperley-Lieb algebra, which is the universal algebra generated by self-adjoint operators $1, U_i$, $i \geq 1$, with relations
\renewcommand{\arraystretch}{1.5}
$$\begin{array}{c}
U_i^2 = \delta U_i, \qquad \qquad
U_i U_j = U_j U_i, \quad |i-j|>1, \\
U_i U_{i+1} U_i - U_i = U_{i+1} U_i U_{i+1} - U_{i+1}, \\
\left( U_i - U_{i+2} U_{i+1} U_i + U_{i+1} \right) \left( U_{i+1} U_{i+2} U_{i+1} - U_{i+1} \right) = 0.
\end{array}$$
\renewcommand{\arraystretch}{1}

For the $A_2$-planar algebra $A_2$-$PTL$, the labels for the internal discs in are now just $A_2$-annular tangles. Let $A_2$-$AnnTL(\sigma,\sigma')$ be the set of all basis $A_2$-annular $(\sigma,\sigma')$-tangles with no other internal discs. Elements of $A_2$-$AnnTL(\sigma,\sigma')$ define elements of $A_2$-$ATL$ by passing to the quotient of $A_2$-$FATL$ by relations K1-K3. The objects of $A_2$-$ATL$ are pairs $(\sigma,\overline{a})$ of a sign string $\sigma$ and a colour $a \in \{ 0,1,2 \}$. The vector space $A_2$-$ATL_{\sigma}$ has as basis the set of $A_2$-annular $\sigma$-tangles with no internal discs, contractible circles, digons or embedded squares. However, non-contractible circles are allowed, making each vector space $A_2$-$ATL((\sigma,\overline{a}),(\sigma',\overline{a'}))$ infinite dimensional. Multiplication in $A_2$-$ATL_{\sigma}(\delta)$ is by composition of tangles (where this makes sense), then reducing the resulting tangle using relations K1-K3 to remove closed loops (K1), digons (K2) or embedded squares (K3).

\begin{figure}[tb]
\begin{center}
  \includegraphics[width=90mm]{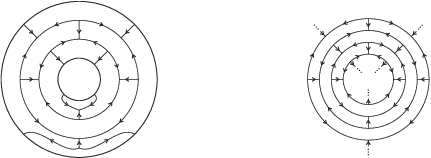}\\
 \caption{A basis $A_2$-annular $--++$-tangle containing hexagons, and the possibility of an infinite number of hexagons}\label{fig:ATL1}
\end{center}
\end{figure}

For $\sigma, \sigma'$ be sign strings of length $\geq 3$, the algebras $A_2$-$ATL((\sigma,\overline{a}),(\sigma,\overline{a'}))$ are also infinite dimensional due to the possibility of an infinite number of embedded hexagons in basis tangles in the annular picture, as illustrated in Figure \ref{fig:ATL1}.

We have a notion of the rank of a tangle. A minimal cut loop $\gamma$ in an annular $(\sigma,\sigma')$-tangle $T$ will be a clockwise closed path which encloses the distinguished internal disc and crosses the least number of strings. We associate a weight $w_{\gamma} = (t_1,t_2)$ to a minimal cut loop $\gamma$, where $t_1$ is the number of strings of $T$ that cross $\gamma$ with orientation from left to right, and $t_2$ the number of strings that have orientation from right to left, as we move along $\gamma$ in a complete clockwise loop.
For a weight $(t_1,t_2)$, let $t_{\mathrm{max}} = \mathrm{max} \{ t_1, t_2 \}$ and $t_{\mathrm{min}} = \mathrm{min} \{ t_1, t_2 \}$.
We write $(t_1', t_2') < (t_1, t_2)$, if $t_1' + t_2' < t_1 + t_2$, and if $t_1' + t_2' = t_1 + t_2$ then $(t_1', t_2') < (t_1, t_2)$ if $2t_{\mathrm{max}}' + t_{\mathrm{min}}' < 2t_{\mathrm{max}} + t_{\mathrm{min}}$. The \textbf{rank} of $T$ is then given by the smallest weight $w_{\gamma}$ associated to a minimal cut loop, such that $w_{\gamma} \leq w_{\gamma'}$ for all other minimal cut loops $\gamma'$.

\begin{figure}[bt]
\begin{center}
  \includegraphics[width=115mm]{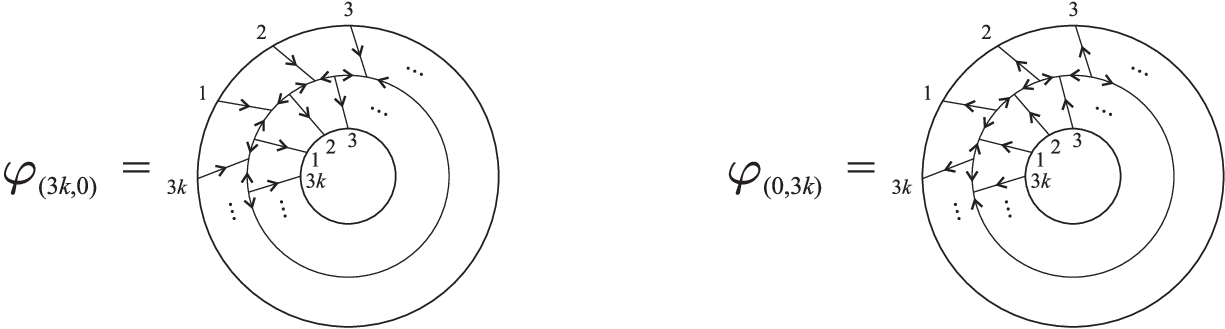}
 \caption{$\varphi_{(3k,0)}$ and $\varphi_{(0,3k)}$} \label{fig:ATL2}
\end{center}
\end{figure}

We define the weight $w(\sigma)$ of a sign string $\sigma$ to be $(s_+,s_-)$, where $s_{\pm}$ is the number of signs `$\pm$' in $\sigma$. The weights $w(\sigma)$ are giving the same ordering as the weights $w_{\gamma}$ of minimal cut loops.
For an $A_2$-planar algebra $P$, we denote by $A_2$-$AP_{\sigma}^{(t_1,t_2)}$ the linear span in the algebra $A_2$-$AP_{\sigma}$ of all labelled $A_2$-annular $\sigma$-tangles with rank $(t_1',t_2') < (t_1,t_2)$, for any rank $(t_1, t_2) \leq w(\sigma)$. Since the rank cannot increase under composition of tangles, $A_2$-$AP_{\sigma}^{(t_1,t_2)}$ is a two-sided ideal in $A_2$-$AP_{\sigma}$.
Note that the quotient of $A_2$-$AP_{\sigma}$ by the ideal $A_2$-$AP_{\sigma}^{(t_1,t_2)}$ is not in general finite dimensional, for $(t_1, t_2) \leq w(\sigma)$. For example, consider the quotient of $A_2$-$AP_{\sigma}$ by $A_2$-$AP_{\sigma}^{(3k,0)}$ (or $A_2$-$AP_{\sigma}^{(0,3k)}$), for $3 \leq 3k \leq w(\sigma)$. The elements $\varphi_{(3k,0)}$, $\varphi_{(0,3k)}$ in Figure \ref{fig:ATL2} have ranks $(3k,0)$, $(0,3k)$ respectively, and can be composed an infinite number of times, but the resulting tangle cannot be reduced using K1-K3.

A $P$-module $V$ is indecomposable if and only if $V_{\sigma}$ is an indecomposable $A_2$-$AP_{\sigma}$-module for each sign string $\sigma$.

The \textbf{weight} $wt(V)$ of a $P$-module $V$ is the smallest weight $w(\sigma)$ for which $V_{\sigma}$ is non-zero. Elements of $V_{\sigma'}$ for any $\sigma'$ such that $w(\sigma') = wt(V)$ will be called lowest weight vectors in $V$, and $V_{\sigma'}$ is an $A_2$-$AP_{\sigma'}$-module which we call a \textbf{lowest weight module}.

Let $+^m$ denote the sign string $++ \cdots +$ ($m$ copies), and $-^n = -- \cdots -$ ($n$ copies). Note that $V_{\sigma} \cong V_{+^m -^n}$ for all sign strings $\sigma$ which are a permutation of the sign string $+^m -^n$, $m \equiv n \textrm{ mod } 3$, where an isomorphism is given by using the braiding to permute the vertices on the outer disc of every element in $V_{\sigma}$ so that the outer disc has pattern $+^m -^n$. In this situation we will say that $\sigma$ is of type $+^m -^n$.

\subsection{Hilbert $P$-modules} \label{Sect:Hilbert_P-modules}

We define Hilbert $P$-modules, parallel to \cite[Section 3]{jones:2001}.
Let $P$ be a \textbf{$C^{\ast}$-$A_2$-planar algebra}, that is, a non-degenerate finite-dimensional $A_2$-planar $\ast$-algebra with positive definite partition function $Z(T): \otimes_k P_{\sigma_k} \longrightarrow P_{\varnothing} = \mathbb{C}$, for any $\varnothing$-tangle $T$.
The $\ast$-algebra structure on $P$ induces a $\ast$-structure on $A_2$-$AP$, where the involution $\ast$ is defined by reflecting an $A_2$-annular $(\sigma,\sigma')$-tangle $T$ about a circle halfway between the inner and outer disc, and reversing the orientation. $T^{\ast}$ will be an $A_2$-annular $(\sigma',\sigma)$-tangle. If $P$ is a $C^{\ast}$-$A_2$-planar algebra this defines an antilinear involution $\ast$ on $A_2$-$FAP$ by taking the $\ast$ of the underlying unlabelled tangle for a labelled tangle $T$, replacing the labels of $T$ by their $\ast$'s, and extending by antilinearity. Since $P$ is an $A_2$-planar $\ast$-algebra, all the $A_2$-planar relations are preserved under $\ast$ on $A_2$-$FAP$, so $\ast$ passes to an antilinear involution on the algebra $A_2$-$AP$. In particular, all the $A_2$-$AP_{\sigma}$ are $\ast$-algebras.

Let $P$ be a $C^{\ast}$-$A_2$-planar algebra. A $P$-module $H$ will be called a \textbf{Hilbert $P$-module} if each $H_{\sigma}$ is a finite dimensional Hilbert space with an invariant inner-product $\langle \cdot ,  \cdot \rangle$, i.e.
\begin{equation} \label{eqn:inner-product_on_Hilbert-module}
\langle av, w \rangle = \langle v, a^{\ast} w \rangle,
\end{equation}
for all $v,w, \in H$ and $a \in A_2$-$AP$.
As in the case of $A_1$-planar algebras, a $P$-submodule of a Hilbert $P$-module is a Hilbert $P$-module. Also, the orthogonal complement of a $P$-submodule is a $P$-module, so that indecomposability and irreducibility are the same for Hilbert $P$-modules. The following Lemma follows immediately from Lemma 3.4 in \cite{jones:2001}:

\begin{Lemma} \label{Lemma:SU(3)Lemma3.4}
Let $P$ be an $A_2$-$C^{\ast}$-planar algebra and $H$ a Hilbert $P$-module. If $W \subseteq H_{\sigma}$ is an irreducible $A_2$-$AP_{\sigma}$-submodule of $H_{\sigma}$ for some $\sigma$, then $A_2$-$AP(W)$ is an irreducible $P$-submodule of $H$.
\end{Lemma}

From the invariance of $\langle \cdot ,  \cdot \rangle$ it is easy to see that if $V$ and $W$ are orthogonal $A_2$-$AP_{\sigma}$ invariant subspaces of $H_{\sigma}$ for some $\sigma$, then $A_2$-$AP(V)$ is orthogonal to $A_2$-$AP(W)$.
As in \cite{jones:2001}, an irreducible Hilbert $P$-module $H$ is determined by its lowest weight modules. In particular $H$ is determined by a lowest weight module $H_{+^m -^n}$.

We now determine which $A_2$-$AP_{\sigma}$-modules can be lowest weight modules.

\begin{Lemma} \label{Lemma:P-modulesA1}
Let $P$ be an $A_2$-$C^{\ast}$-planar algebra and $H$ a Hilbert $P$-module of rank $(t_1,t_2)$. For $\sigma$ such that $w(\sigma) = wt(H)$, any element $w \in H_{\sigma}$ can be written, up to a scalar, as $aw$ for some $a \in A_2\textrm{-}AP_{\sigma}$ with $\mathrm{rank}(a) = \mathrm{rank}(w)$.
\end{Lemma}
\emph{Proof:}
First form $w w^{\ast} \in A_2\textrm{-}AP_{\sigma}$. Then dividing out by the relations K1-K3 we obtain a linear combination of elements in $A_2$-$AP_{\sigma}$, and we remove those elements that have rank $< (t_1,t_2)$. Ignoring the scalar factor we are left with a single element $a \in A_2\textrm{-}AP_{\sigma}$ with $\mathrm{rank}(a) = (t_1,t_2)$. If we form $aw$, then dividing out by K1-K3 we obtain $aw = \mu w + \sum_i \mu_i w_i$, where $\mu,\mu_i \in \mathbb{C}$ and $w_i \in H_{\sigma}$ with $\mathrm{rank}(w_i) < (t_1,t_2)$ for each $i$. Then in $H$ the $w_i$ are all zero, so that $\mu^{-1}aw=w$.
\hfill
$\Box$

Then the proof of Lemma 3.8 in \cite{jones:2001} gives the following result:

\begin{Lemma} \label{Lemma:SU(3)Lemma3.8}
Let $P$ be an $A_2$-$C^{\ast}$-planar algebra and $H$ a Hilbert $P$-module. For $\sigma$ any sign string of type $+^{t_1} -^{t_2}$, let $H^{(t_1,t_2)}_{\sigma}$ be the $A_2$-$AP_{\sigma}$-submodule of $H_{\sigma}$ spanned by the $\sigma$-graded pieces of all $P$-submodules with rank $< (t_1,t_2)$. Then
$$(H^{(t_1,t_2)}_{\sigma})^{\bot} = \bigcap_{a \in A_2\textrm{-}AP^{(t_1,t_2)}_{\sigma}} \mathrm{ker}(a).$$
\end{Lemma}

Thus we see that the lowest weight modules of an irreducible $P$-module of rank $(t_1,t_2)$ are $A_2\textrm{-}AP_{\sigma}/A_2\textrm{-}AP^{(t_1,t_2)}_{\sigma}$-modules, for all $\sigma$ of type $+^{t_1} -^{t_2}$. In fact, for an $A_2$-$C^{\ast}$-planar algebra $P$, to determine all Hilbert $P$-modules it is sufficient to first determine the algebras $A_2\textrm{-}AP_{+^{t_1} -^{t_2}}/A_2\textrm{-}AP^{(t_1,t_2)}_{+^{t_1} -^{t_2}}$ and their irreducible modules, where $A_2\textrm{-}AP^{\overline{0},(t_1,t_2)}_{\sigma}$ is the two-sided ideal in $A_2\textrm{-}AP_{\sigma}^{\overline{0}}$ given by the linear span of all labeled $A_2$-annular $\sigma$-tangles in $A_2\textrm{-}AP_{\sigma}^{\overline{0}}$ with rank $< (t_1,t_2)$.

\subsection{Irreducible $A_2$-$PTL$-modules} \label{sect:Irr_PTL-modules}

In this section we will describe all zero-weight irreducible modules, parallel to the construction of zero-weight modules in \cite{jones:2001}.

\begin{figure}[htb]
\begin{center}
  \includegraphics[width=100mm]{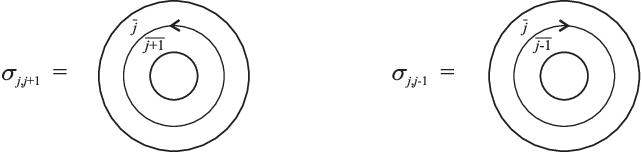}
 \caption{$\sigma_{j,j+1}$ and $\sigma_{j,j-1}$} \label{fig:ATL4}
\end{center}
\end{figure}

For $i \in \{ 0,1,2 \}$, $\varepsilon \in \{ \pm \}$, and non-negative integers $k_1 \equiv k_2 \textrm{ mod } 3$, let
$$\sigma^{(i,\varepsilon)}_{(k_1,k_2)} = (\sigma_{i,i\varepsilon1} \sigma_{i\varepsilon1,i\varepsilon2} \cdots \sigma_{i\varepsilon k_1-\varepsilon1,i\varepsilon k_1})(\sigma_{i\varepsilon k_1,i\varepsilon k_1-\varepsilon1} \sigma_{i\varepsilon k_1-\varepsilon1,i\varepsilon k_1-\varepsilon2} \cdots \sigma_{i\varepsilon k_1-\varepsilon k_2\varepsilon1,i\varepsilon k_1-\varepsilon k_2}),$$
where the $0$-tangles $\sigma_{j,j \pm 1}$ are illustrated in Figure \ref{fig:ATL4}, $j \in \{ 0,1,2 \}$.
The algebra $A_2\textrm{-}ATL_{\varnothing,\overline{a}}$ is generated by the $\varnothing$-tangles $\sigma^{(i,\pm)}_{(k_1,k_2)}$.
Let $H$ be an irreducible Hilbert $A_2$-$PTL$-module of lowest weight zero. For any $i \in \{ 0,1,2 \}$, $\varepsilon \in \{ \pm \}$, and non-negative integers $k_1 \equiv k_2 \textrm{ mod } 3$, $Z(\sigma^{(i,\pm)}_{(k_1,k_2)}) = \lambda \in \mathbb{C}$ since $Z$ is a positive definite partition function for any $\varnothing$-tangle, and $Z(\sigma^{(i,\pm)\ast}_{(k_1,k_2)}) = Z(\sigma^{(i,\pm)}_{(k_1,k_2)})^{\ast} = \overline{\lambda}$. Then since the partition function is multiplicative on disconnected components, we see that
$Z((\sigma^{(i,\pm)}_{(k_1,k_2)}-\lambda \mathbf{1}_{\varnothing})^{\ast} (\sigma^{(i,\pm)}_{(k_1,k_2)}-\lambda \mathbf{1}_{\varnothing})) = Z(\sigma^{(i,\pm)\ast}_{(k_1,k_2)})Z(\sigma^{(i,\pm)}_{(k_1,k_2)}) - \overline{\lambda} Z(\sigma^{(i,\pm)}_{(k_1,k_2)}) - \lambda Z(\sigma^{(i,\pm)\ast}_{(k_1,k_2)}) + |\lambda|^2 Z(\mathbf{1}_{\varnothing}) = 0$,
thus $\sigma^{(i,\pm)}_{(k_1,k_2)}$ acts by the scalar $\lambda$ in any Hilbert $A_2$-$PTL$-module of lowest weight zero.

Suppose $\sigma^{(0,+)}_{(3,0)}$, $\sigma^{(0,+)}_{(1,1)}$, $\sigma^{(0,-)}_{(1,1)}$ act as scalars $a,b_+,b_- \in \mathbb{C}$ respectively. Since $\sigma^{(0,+)}_{(0,3)} = \sigma^{(0,+)\ast}_{(3,0)}$, we see that $\sigma^{(0,+)}_{(0,3)}$ acts as the scalar $\overline{a}$, and since $\sigma^{(0,\pm)}_{(1,1)}$ are self-adjoint we must have $b_{\pm} \in \mathbb{R}$. Consider the $\varnothing$-tangle $\sigma^{(0,+)}_{(0,3)}\sigma^{(0,+)\ast}_{(3,0)}$ which acts as the scalar $|a|^2$. However, this tangle contains three copies of the tangle $\sigma^{(0,+)}_{(1,1)}$, which acts as a scalar $b_+^3$. Thus $b_+$ must be the real cubic root of $|a|^2$. If we let $a = \beta^3$ for some $\beta \in \mathbb{C}$, then $b_+ = |\beta|^2$. Similarly, from considering the $\varnothing$-tangle $\sigma^{(0,+)\ast}_{(0,3)}\sigma^{(0,+)}_{(3,0)}$, we find that $b_- = |\beta|^2 = b_+$. The scalar $\beta$ is unique up to a choice of third root of unity, but the Hilbert $A_2$-$PTL$-module $H$ does not depend on this choice as the only factors of $\beta$ which appear are of the form $\beta^3$ or $|\beta|^2$.
Then we obtain the following result:

\begin{Prop} \label{Thm:SU(3)5.12}
An irreducible Hilbert $A_2$-$PTL$-module $H$ of weight zero in which the maps $\sigma^{(i,+)}_{(k_1,k_2)}$, $\sigma^{(i,-)}_{(k_1,k_2)}$, $i \in \{ 0,1,2 \}$, act as the complex number $\beta^{k_1} \overline{\beta}^{k_2}$, $\beta^{k_2} \overline{\beta}^{k_1}$ respectively, for some fixed $\beta \in \mathbb{C}$, is determined up to isomorphism by the dimensions of $H_{\varnothing,\overline{a}}$, $a \in \{ 0,1,2 \}$, and the number $\beta$, where we require $|\beta| \leq \alpha$.
\end{Prop}
\emph{Proof:}
The uniqueness of the $A_2$-$PTL$-module is a consequence of the fact that an irreducible Hilbert $P$-module $H$ is determined by its lowest weight modules (see Section \ref{Sect:Hilbert_P-modules}), since at least one of $H_{\varnothing,\overline{0}}$, $H_{\varnothing,\overline{1}}$ and $H_{\varnothing,\overline{2}}$ is non-zero. Let $E_1$, $E_2$ be the tangles
\begin{center}
  \includegraphics[width=110mm]{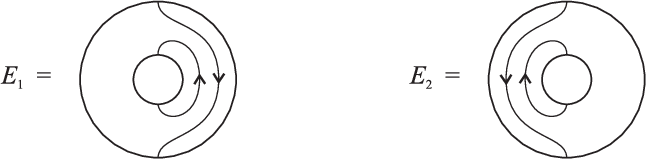}
\end{center}
so that $\alpha^{-1} E_1$, $\alpha^{-1} E_2$ are projections. Then since $E_1 E_2 E_1 = |\beta|^2 E_1$ we have $||\alpha^{-1}E_1 \cdot \alpha^{-1}E_2 \cdot \alpha^{-1}E_1|| = |\beta|^2 \alpha^{-2} ||\alpha^{-1}E_1||$ so that $1 \geq |\beta|^2 \alpha^{-2}$. Hence $|\beta| \leq \alpha$.
\hfill
$\Box$

For $\beta = \alpha$, $V_{\sigma}^{\alpha} = A_2\textrm{-}PTL_{\sigma}$ (since when $\beta = \alpha$ there is no distinction between contractible and non-contractible circles). For $\alpha > 3$ (which corresponds to $\delta > 2$), the inner product is positive definite by \cite[Lemma 3.10]{evans/pugh:2009iii} and \cite[Theorem 3.6]{wenzl:1988}, and $H_{\sigma}^{\alpha} = V_{\sigma}^{\alpha}$ is a Hilbert $A_2$-$PTL$-module. For $0 < \alpha \leq 3$, if the inner product is positive semi-definite on $V_{\sigma}^{\alpha}$ we let $H_{\sigma}^{\alpha}$ be the quotient of $V_{\sigma}^{\alpha}$ by the subspace of vectors of length zero; otherwise $H_{\sigma}^{\alpha}$ does not exist.

Now consider the case when $0 < |\beta| < \alpha$. We define the set $Th_{\sigma}$ to be the set of all $(\sigma,\overline{0})$-tangles with no contractible circles and at most two non-contractible circles. For each $\beta$ we form the graded vector space $V^{\beta}$, where $V^{\beta}_{\sigma}$ has basis $Th_{\sigma}$, and we equip it with an $A_2$-$PTL$-module structure as follows. Let $T \in A_2\textrm{-}ATL(\sigma',\sigma)$ and $R \in A_2\textrm{-}ATL_{\sigma}$. We from the tangle $TR$ and reduce it using K1-K3, so that $TR = \sum_j \delta^{b_j} \alpha^{c_j} (TR)_j$, for some basis $A_2$-annular $(\sigma',\overline{0})$-tangles $(TR)_j$, where $b_j$, $c_j$ are non-negative integers. Let $\sharp^a_j$, $\sharp^c_j$ denote the number of non-contractible circles in the tangle $(TR)_j$ which have anti-clockwise, clockwise orientation respectively. We define integers $d_j$, $f_j$ and $g_j$ as follows: $d_j = \textrm{min}(\sharp^a_j, \sharp^c_j)$, $f_j = \sharp^a_j - \sharp^c_j - \gamma_{f_j}$ if $\sharp^a_j \geq \sharp^c_j$ and $f_j = 0$ otherwise, and $g_j = \sharp^c_j - \sharp^a_j - \gamma_{g_j}$ if $\sharp^a_j \leq \sharp^c_j$ and $g_j = 0$ otherwise, where $\gamma_{f_j}, \gamma_{g_j} \in \{ 0,1,2 \}$ such that $f_j, g_j \equiv 0 \textrm{ mod }3$. Then we set $T(R) = \sum_j \delta^{b_j} \alpha^{c_j} \beta^{d_j + f_j} \overline{\beta}^{d_j + g_j} (\widehat{TR})_j$, where $(\widehat{TR})_j$ is the tangle $(TR)_j$ with $d_j + f_j$ anti-clockwise non-contractible circles removed, and $d_j + g_j$ clockwise ones removed.

\begin{Prop}
The above definition make $V^{\beta}$ into an $A_2$-$PTL$-module of weight zero in which $\sigma^{(a)}_{(k_1,k_2)} = \beta^{k_1} \overline{\beta}^{k_2}$ for $a=0,1,2$.
\end{Prop}

The choice of $(\sigma,\overline{0})$-tangles rather than $(\sigma,\overline{1})$- or $(\sigma,\overline{2})$-tangles to define $V^{\beta}$ was arbitrary. For these other two choices, the maps $T \rightarrow \beta^{-1} T \sigma_{01}$, $T \rightarrow \overline{\beta}^{-1} T \sigma_{02}$ respectively would define isomorphisms from those modules to the one defined above.

\begin{Def} \label{Def:SU(3)5.17}
Given two tangles $S,T \in Th_{\sigma}$, we reduce $T^{\ast} S$ using K1-K3 so that $T^{\ast} S = \sum_j \delta^{b_j} \alpha^{c_j} (T^{\ast}S)_j$, where $(T^{\ast}S)_j$ are basis tangles in $A_2$-$ATL_{\varnothing}$. Define $d_j$, $f_j$ and $g_j$ for each $(T^{\ast}S)_j$ as above. We define an inner-product by $\langle S,T \rangle = \sum_j \delta^{b_j} \alpha^{c_j} \beta^{d_j + f_j} \overline{\beta}^{d_j + g_j}$.
\end{Def}

Invariance of this inner-product follows from the fact that $T^{\ast}S = \langle S,T \rangle \mathbf{1}_{\varnothing}$ where $\mathbf{1}_{\varnothing}$ is the annular $(\overline{0}:\overline{0})$-tangle with no strings at all. When the above inner-product is positive semi-definite, we define the Hilbert $A_2$-$PTL$-module $H^{\beta}$ of weight zero to be the quotient of $V^{\beta}$ by the subspace of vectors of length zero. Otherwise $H^{\beta}$ does not exist.

\begin{Prop}
For the above Hilbert $A_2$-$PTL$-module $H^{\beta}$ of weight zero, the dimension of $H^{\beta}_{\varnothing,\overline{a}}$ is either 0 or 1 for any $\beta \in \mathbb{C} \setminus \{ 0 \}$.
\end{Prop}
\emph{Proof:}
For $a=0$ the result is trivial since $V^{\beta}_{\varnothing,\overline{0}}$ is the linear span of the empty tangle $\mathbf{1}_{\varnothing}$ given in Defn. \ref{Def:SU(3)5.17}. For $a=1$, $V^{\beta}_{\varnothing,\overline{1}} = \mathrm{span}(\sigma_{10}, \sigma_{12} \sigma_{20})$. Let $w = |\beta|^2 \sigma_{12} \sigma_{20} - \beta^3 \sigma_{10}$. Then
\begin{eqnarray*}
\langle w,w \rangle & = & |\beta|^4 \langle \sigma_{12} \sigma_{20}, \sigma_{12} \sigma_{20} \rangle - |\beta|^2 \overline{\beta}^3 \langle \sigma_{12} \sigma_{20}, \sigma_{10} \rangle - \beta^3 |\beta|^2 \langle \sigma_{10}, \sigma_{12} \sigma_{20} \rangle + |\beta|^6 \langle \sigma_{10}, \sigma_{10} \rangle \\
& = & |\beta|^4(|\beta|^4) - |\beta|^2 \overline{\beta}^3 \beta^3 - \beta^3 |\beta|^2 \overline{\beta}^3 + |\beta|^6 |\beta|^2 \;\; = \;\; 0.
\end{eqnarray*}
Then $\sigma_{10} = |\beta|^2 \beta^{-3} \sigma_{12}\sigma_{20} = \overline{\beta} \beta^{-2} \sigma_{12}\sigma_{20}$ in $H^{\beta}_{\varnothing,\overline{1}}$. Similarly when $a=2$, $\sigma_{21} \sigma{10} = \overline{\beta}^2 \beta^{-1} \sigma_{20}$.
\hfill
$\Box$

So we may write a basis for $H^{\beta}$ which does not contain any clockwise non-contractible circles, where we replace every  $\sigma_{10}$ by $\overline{\beta} \beta^{-2} \sigma_{12}\sigma_{20}$ and every $\sigma_{21} \sigma_{10}$ by $\overline{\beta}^2 \beta^{-1} \sigma_{20}$.

\begin{Prop}
The Hilbert $A_2$-$PTL$-module $H^{\beta}$, $|\beta| < \alpha$, is irreducible.
\end{Prop}
\emph{Proof:}
Since $H^{\beta}_{\varnothing,\overline{a}}$ is at most one-dimensional it must be irreducible, for each $a \in \{ 0,1,2 \}$. The maps $\sigma_{j,j+1}$ moves a non-zero element in $H^{\beta}_{\varnothing,\overline{j}}$ to an element in $H^{\beta}_{\varnothing,\overline{j+1}}$, and hence the lowest weight module $H^{\beta}_{\varnothing} = H^{\beta}_{\varnothing,\overline{0}} \oplus H^{\beta}_{\varnothing,\overline{1}} \oplus H^{\beta}_{\varnothing,\overline{2}}$ is irreducible as an $A_2$-$ATL_{\varnothing}$-module. Since $H^{\beta} = A_2\textrm{-}ATL(H^{\beta}_{\varnothing})$, the result follows from Lemma \ref{Lemma:SU(3)Lemma3.4}.
\hfill
$\Box$

Now we consider the case when $\beta = 0$. For each $\sigma$, the set $Th^{\overline{a}}_{\sigma}$ is defined to be the set of all $(\sigma,\overline{a})$-tangles with no contractible or non-contractible circles at all. The cardinality of $Th^{\overline{a}}_{\overline{b}}$ is $\delta_{a,b}$. We form the graded vector space $V^{0,\overline{a}}$, where $V^{0,\overline{a}}_{\sigma}$ has basis $Th^{\overline{a}}_{\sigma}$. We equip it with an $A_2$-$PTL$-module structure of lowest weight zero as follows. Let $T \in A_2\textrm{-}ATL(\sigma',\sigma)$ and $R \in Th^{\overline{a}}_{\sigma}$. We form $TR$ and reduce it using K1-K3, so that $TR = \sum_j \delta^{b_j} \alpha^{c_j} (TR)_j$ as in the case $0 < |\beta| < \alpha$. We define $(\widehat{TR})_j$ to be zero if there are any non-contractible circles in $(TR)_j$, and $(TR)_j$ otherwise. Then $T(R) = \sum_j \delta^{b_j} \alpha^{c_j} (\widehat{TR})_j$.

\begin{Prop}
The above definition make $V^{0,\overline{a}}$ into an $A_2$-$PTL$-module of weight zero in which $\sigma_{j,j \pm 1} = 0$ for $j = 0,1,2 \textrm{ mod } 3$.
\end{Prop}

\begin{Def} \label{Def:SU(3)5.22}
Given basis tangles $S,T \in Th^{\overline{a}}_{\sigma}$, we reduce $T^{\ast} S$ using K1-K3 so that $T^{\ast} S = \sum_j \delta^{b_j} \alpha^{c_j} (T^{\ast}S)_j$ for basis $(\overline{a}:\overline{a})$-tangles $(T^{\ast}S)_j$. We define $\langle S,T \rangle_j$ to be 0 if there are any non-contractible circles in $(T^{\ast}S)_j$, and 1 otherwise. Then we define an inner-product by $\langle S,T \rangle = \sum_j \delta^{b_j} \alpha^{c_j} \langle S,T \rangle_j$.
\end{Def}

This inner-product is invariant as in the case $0 < |\beta| < \alpha$. Again, if the inner product is positive semi-definite we define $H^{0,\overline{a}}$ to be the quotient of $V^{0,\overline{a}}$ by the subspace of vectors with length zero; otherwise $H^{0,\overline{a}}$ does not exist.
The Hilbert $A_2$-$PTL$-module $H^{0,\overline{a}}$, $a \in \{ 0,1,2 \}$, is irreducible,
where the proof is as for $H^{\beta}$.

\section{The $A_2$-graph planar algebra of an $\mathcal{ADE}$ graph} \label{Sect:A2-GPA}

We now describe the construction of the $A_2$-graph planar algebra $P^{\mathcal{G}}$ for an $\mathcal{ADE}$ graph. We will then determine a partial decomposition of $P^{\mathcal{G}}$ into irreducible $A_2$-$PTL$-modules.

Let $\mathcal{G}$ be any finite $SU(3)$ $\mathcal{ADE}$ graph with vertex set $\mathfrak{V}^{\mathcal{G}} = \mathfrak{V}^{\mathcal{G}}_0 \cup \mathfrak{V}^{\mathcal{G}}_1 \cup \mathfrak{V}^{\mathcal{G}}_2$, where $\mathfrak{V}^{\mathcal{G}}_a$ is the set of $a$-coloured vertices of $\mathcal{G}$, $a = 0,1,2$. Let $s_a = |\mathfrak{V}^{\mathcal{G}}_a|$ denote the number of $a$-coloured vertices and let $s = |\mathfrak{V}^{\mathcal{G}}|$, the total number of vertices of $\mathcal{G}$. For a three-colourable graph, $s = s_0 + s_1 + s_2$, and we have $s_1 = s_2$ due to the conjugation property of the $SU(3)$ $\mathcal{ADE}$ graphs. For the non-three-colourable graphs we discard any colouring of the vertices, so that $s = s_a$, $a = 0,1,2$.
Let $\alpha = [3]_q$, $q=e^{i \pi/n}$, be the Perron-Frobenius eigenvalue of $\mathcal{G}$ and let $\phi = (\phi_v)$ be the corresponding eigenvector.

\begin{figure}[tb]
\begin{center}
\includegraphics[width=100mm]{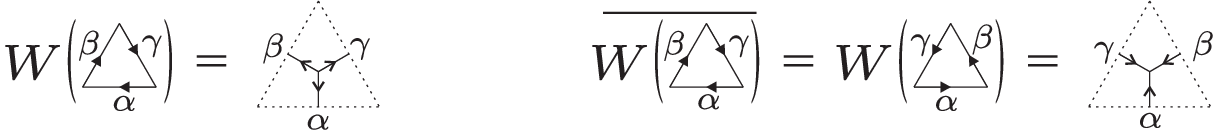}\\
 \caption{Cells associated to trivalent vertices} \label{fig:Oc-Kup}
\end{center}
\end{figure}

Ocneanu \cite{ocneanu:2000ii} defined a cell system $W$ on $\mathcal{G}$, associating a complex number $W \left( \triangle^{(\alpha \beta \gamma)} \right)$, called an Ocneanu cell, to each closed loop of length three $\triangle^{(\alpha \beta \gamma)}$ in $\mathcal{G}$ as in Figure \ref{fig:Oc-Kup}, where $\alpha$, $\beta$, $\gamma$ are edges on $\mathcal{G}$, and such that these cells satisfy two properties, called Ocneanu's type I, II equations respectively, which are obtained by evaluating the Kuperberg relations K2, K3 respectively, using the identification in Figure \ref{fig:Oc-Kup}: \\
$(i)$ for any type I frame \includegraphics[width=16mm]{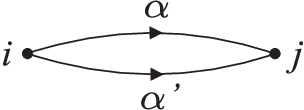} in $\mathcal{G}$ we have \\
\begin{flushright}
\begin{minipage}[b]{14cm}
 \begin{minipage}[t]{9cm}
  \mbox{} \\
 \includegraphics[width=80mm]{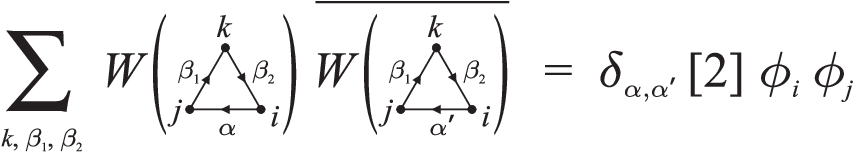}
 \end{minipage}
 \hfill
 \begin{minipage}[t]{1.5cm}
  \mbox{} \\
  \hfill
  \parbox[t]{7mm}{\begin{eqnarray}\label{eqn:typeI_frame}\end{eqnarray}}
 \end{minipage}
\end{minipage}
\end{flushright}
$(ii)$ for any type II frame \includegraphics[width=30mm]{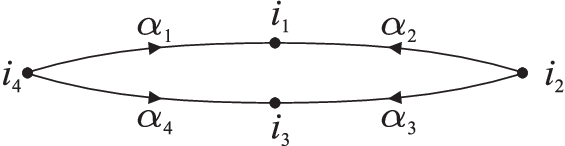} in $\mathcal{G}$ we have \\
\begin{flushright}
\begin{minipage}[b]{14cm}
 \begin{minipage}[b]{9cm}
  \mbox{} \\
 \includegraphics[width=90mm]{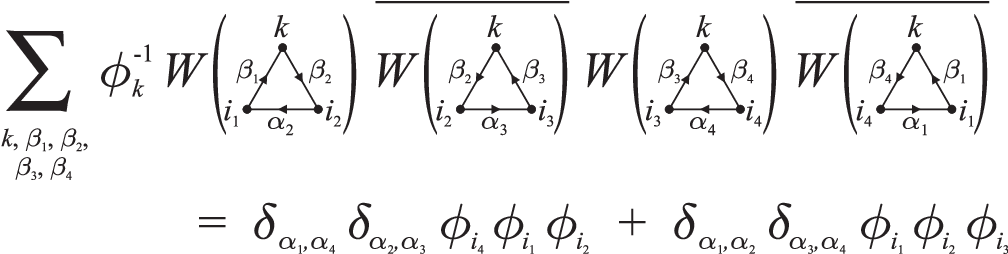}
 \end{minipage}
 \hfill
 \begin{minipage}[b]{1.5cm}
  \mbox{} \\
  \hfill
  \parbox[b]{7mm}{\begin{eqnarray}\label{eqn:typeII_frame}\end{eqnarray}}
 \end{minipage}
\end{minipage}
\end{flushright}
The existence of these cells for the finite $\mathcal{ADE}$ graphs was shown in \cite{evans/pugh:2009i} with the exception of the graph $\mathcal{E}_4^{(12)}$.

We denote by $\mathcal{G}^{\mathrm{op}}$ the reverse graph of $\mathcal{G}$, which is the graph obtained by reversing the direction of every edge of $\mathcal{G}$. For an edge $\gamma$ on $\mathcal{G}$ we will denote the reverse edge on $\mathcal{G}^{\mathrm{op}}$ by $\widetilde{\gamma}$.
Let $P^{\mathcal{G}}_{\sigma} = \bigoplus_{a=0}^{2} P^{\mathcal{G},\overline{a}}_{\sigma}$ where $P^{\mathcal{G},\overline{a}}_{\sigma}$ is the space of closed paths on $\mathcal{G}$, $\mathcal{G}^{\mathrm{op}}$ which begin at vertices in $\mathfrak{V}^{\mathcal{G}}_a$, and with pattern $\sigma$ where a `$-$' denotes that an edge is on $\mathcal{G}$ and `$+$' denotes that an edge is on $\mathcal{G}^{\mathrm{op}}$.
Note that for a non-three-colourable graph $\mathcal{G}$, $P^{\mathcal{G}}$ contains three copies of the the space of closed paths on $\mathcal{G}$, $\mathcal{G}^{\mathrm{op}}$, one for each $a = 0,1,2$.

We now define a presenting map $Z:\mathcal{P}(P^{\mathcal{G}}) \rightarrow P^{\mathcal{G}}_{i,j}$ in a similar way to the definition of $Z:\mathcal{P}(P) \rightarrow P$ in \cite{evans/pugh:2009iii} for a subfactor planar algebra $P$ where there is a flat connection \cite{ocneanu:1988, ocneanu:1991} defined on the graph $\mathcal{G}$. The difference is that in the present definition for $P^{\mathcal{G}}$ we do not need to use the connection, and hence we can define $P^{\mathcal{G}}$ for any $SU(3)$ $\mathcal{ADE}$ graph $\mathcal{G}$ and not just those for which a flat connection exists.
We define a $\ast$-operation on $P^{\mathcal{G}}$ by $\gamma^{\ast} = \widetilde{\gamma} \in P^{\mathcal{G}}_{\sigma^{\ast}}$ for $\gamma \in P^{\mathcal{G}}_{\sigma}$.

\begin{figure}[tb]
\begin{center}
   \includegraphics[width=70mm]{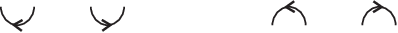}
\caption{left and right cups; left and right caps} \label{fig:cups&caps}
\end{center}
\end{figure}

\begin{figure}[tb]
\begin{center}
   \includegraphics[width=70mm]{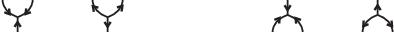}
\caption{incoming and outgoing Y-forks; incoming and outgoing inverted Y-forks} \label{fig:Y-forks}
\end{center}
\end{figure}

Let $T$ be a labelled tangle in $\mathcal{P}^{\mathcal{G}}_{\sigma}$ with $m$ internal discs $D_k$ with pattern $\sigma_k$ and labels $x_k \in P^{\mathcal{G}}_{\sigma_k}$, $k=1, \ldots, m$. We define $Z(T)$ as follows. First, convert all the discs $D_k$ to rectangles (including the outer disc) so that its edges are parallel to the $x,y$-axes, and such that all the vertices on its boundary lie along the top edge of the rectangle.
Next, isotope the strings of $T$ so that each horizontal strip only contains one of the following elements: a rectangle with label $x_k$, a cup, a cap, a Y-fork, or an inverted Y-fork
(see Figures \ref{fig:cups&caps} and \ref{fig:Y-forks}).
For a tangle $T \in \mathcal{P}^{\mathcal{G}}_{\sigma}$ with $l$ horizontal strips $s_l$, where $s_1$ is the highest strip, $s_2$ the strip immediately below it, and so on, we define $Z(T)=Z(s_1) Z(s_2) \cdots Z(s_l)$, which will be an element of $P^{\mathcal{G}}_{\sigma}$.
This algebra is normalized in the sense that for the empty tangle $\bigcirc$, $Z(\bigcirc)=1$.
We will need to show that this definition only depends on $T$, and not on the decomposition of $T$ into horizontal strips.

For any horizontal strip $s$ we have sign strings $\sigma_1$, $\sigma_2$ given by the endpoints of the strings along the top, bottom edge respectively of the strip (we will call these endpoints vertices), where, along the top edge `$+$' is given by a sink and `$-$' by a source, and along the bottom edge `$+$' is given by a source and `$-$' by a sink.
Each vertex along the top (or bottom) with downwards, upwards orientation respectively, corresponds to an edge on $\mathcal{G}$, $\mathcal{G}^{\mathrm{op}}$ respectively. Then the top, bottom edge of the strip is labelled by elements in $P^{\mathcal{G}}_{\sigma_1}$, $P^{\mathcal{G}}_{\sigma_2}$ respectively. Then $Z(s)$ defines an operator $M_s \in \textrm{End}(P^{\mathcal{G}}_{\sigma_2},P^{\mathcal{G}}_{\sigma_1})$.

For a strip $s$ containing one of the elements in Figures \ref{fig:cups&caps} and \ref{fig:Y-forks}, the definition of $Z(s)$ is equivalent to that given in \cite{evans/pugh:2009iii} for a subfactor $A_2$-planar algebra, which we summarize here.
We define annihilation operators $c_l$, $c_r$ by:
\begin{equation} \label{eqn:annihilation-lr}
c_l(\alpha\overline{\beta}) = \delta_{s(\alpha),s(\beta)} \frac{\sqrt{\phi_{r(\alpha)}}}{\sqrt{\phi_{s(\alpha)}}} \, s(\alpha), \qquad c_r(\overline{\beta}\alpha) = \delta_{r(\alpha),r(\beta)} \frac{\sqrt{\phi_{s(\alpha)}}}{\sqrt{\phi_{r(\alpha)}}} \, r(\alpha),
\end{equation}
and creation operators $c_l^{\ast}$, $c_r^{\ast}$ as their adjoints, where $\alpha$ is an edge on $\mathcal{G}$ and $\overline{\beta}$ an edge on $\mathcal{G}^{\mathrm{op}}$, and $(\phi_v)_v$ is the Perron-Frobenius eigenvector for the Perron-Frobenius eigenvalue $\alpha$ of $\mathcal{G}$.
We normalize $(\phi_v)$ so that $\sum_{v \in \mathfrak{V}^{\mathcal{G}}_0} \phi_v^2 = 1$.
Define the following fork operators $\curlyvee$, $\overline{\curlyvee}$, $\curlywedge$ and $\overline{\curlywedge}$ by:
\begin{eqnarray}
\curlyvee(\overline{\alpha}) & = & \frac{1}{\sqrt{\phi_{s(\alpha)} \phi_{r(\alpha)}}} \sum_{\beta_1,\beta_2} W(\triangle_{ijk}^{(\alpha,\beta_1,\beta_2)}) \beta_1 \beta_2, \label{eqn:Yfork(in)} \\
\overline{\curlyvee}(\alpha) & = & \frac{1}{\sqrt{\phi_{s(\alpha)} \phi_{r(\alpha)}}} \sum_{\beta_1,\beta_2} \overline{W(\triangle_{ijk}^{(\alpha,\beta_1,\beta_2)})} \overline{\beta_1} \overline{\beta_2}, \label{eqn:Yfork(out)}
\end{eqnarray}
and $\curlywedge = \overline{\curlyvee}^{\ast}$ and $\overline{\curlywedge} = \curlyvee^{\ast}$.
Then $Z$ assigns to the left, right caps the annihilation operators $c_l$, $c_r$ respectively, given by (\ref{eqn:annihilation-lr}), to the left, right cups the creation operators $c_l^{\ast}$, $c_r^{\ast}$ respectively, to the incoming, outgoing Y-forks the operators $\curlyvee$, $\overline{\curlyvee}$ respectively given in (\ref{eqn:Yfork(in)}), (\ref{eqn:Yfork(out)}), and to the incoming, outgoing inverted Y-forks the operators $\curlywedge$, $\overline{\curlywedge}$ respectively.

For a strip $s$ containing a rectangle with label $x = \sum_{\gamma} \lambda_{\gamma} \gamma$, where $\gamma$ are single paths in $P^{\mathcal{G}}_{\sigma}$, we define the operator $M_{s} = Z(s)$ differently to the definition of $Z(s)$ given in \cite{evans/pugh:2009iii} for a subfactor $A_2$-planar algebra. The definition here is simpler as it does not involve the connection. Let $p$, $p'$ be the number of vertical strings to the left, right respectively of the rectangle in strip $s$, with orientations given by the sign strings $\sigma^{(p)}$, $\sigma^{(p')}$ respectively.
Then $\sum_{\gamma, \mu, \nu} \lambda_{\gamma} (\mu \cdot \gamma \cdot \nu, \mu \cdot \nu)$ defines a matrix $M_{s}$, where the summation is over all paths $\mu$ of length $p$ with edges on $\mathcal{G}$, $\mathcal{G}^{\mathrm{op}}$ as dictated by the sign string $\sigma^{(p)}$, and paths $\nu$ of length $p'$ with edges on $\mathcal{G}$, $\mathcal{G}^{\mathrm{op}}$ as dictated by the sign string $\sigma^{(p')}$.

\begin{Thm}
The above definition of $Z(T)$ for any $A_2$-planar tangle $T$ makes $P^{\mathcal{G}} = \bigcup_{\sigma} P^{\mathcal{G}}_{\sigma}$ into an $A_2$-$C^{\ast}$-planar algebra, the $A_2$-graph planar algebra,
with $\mathrm{dim}(P^{\mathcal{G},\overline{a}}_{\varnothing}) = n_a$, $a=0,1,2$, and parameter $[3]$.
\end{Thm}
\emph{Proof:}
This follows as in the proof of Theorem 5.4 in \cite{evans/pugh:2009iii}, where the only difference occurs for isotopies of the tangle which involve rectangles. Here the invariance of $Z$ under these isotopies is simpler as the connection is not used -- it follows from the fact that $Z$ assigns vertices to each region in the tangle $T$, and any such assignment does not change under these isotopies.
\hfill
$\Box$

\begin{figure}[bt]
\begin{center}
  \includegraphics[width=35mm]{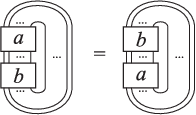}\\
 \caption{$\mathrm{tr}(ab) = \mathrm{tr}(ba)$}\label{fig:tr(ab)=tr(ba)}
\end{center}
\end{figure}

The partition function $Z:\mathcal{P}^{\mathcal{G},\overline{a}}_{\varnothing} \longrightarrow \mathbb{C}$ is defined as the linear map which takes the basis path $v \in \mathfrak{V}^{\mathcal{G}}_0$ to $\phi_v^2$, i.e. the map whose matrix entries, labelled by the vertices of $\mathcal{G}$, are $\phi_v^2$. Thus there is a multiplicative factor $\phi_v^2$ for the external region, which is required for spherical isotopy.
As in the case of the planar algebra of a bipartite graph \cite{jones:2000}, the partition function of a closed labelled tangle $T$ depends only on $T$ up to isotopies of the 2-sphere.
The partition function of an empty tangle is equal to $\sum_{v \in \mathfrak{V}^{\mathcal{G}}_0} \phi_v^2 = 1$.
Then we have the analogue of Theorem 3.6 \cite{jones:2000} for the $A_2$-planar algebra of an $\mathcal{ADE}$ graph:

\begin{Thm}
Let $P^{\mathcal{G}}$ be the $A_2$-graph planar algebra of an $\mathcal{ADE}$ graph $\mathcal{G}$, with (normalized) Perron-Frobenius eigenvector $(\phi_v)$. Then for $x \in P^{\mathcal{G}}_{\sigma\sigma^{\ast}}$, $\mathrm{tr}(x) = [3]^{-|\sigma|} Z(\widehat{x})$ defines a normalized trace on the union of the $P$'s, where $\widehat{x}$ is any $\varnothing$-tangle obtained from $x$ by connecting the first $|\sigma|$ boundary points to the last $|\sigma|$. The scalar product $\langle x, y \rangle = \mathrm{tr} (x^{\ast} y)$ is positive definite.
\end{Thm}
\emph{Proof:}
The normalization makes the definition of the trace consistent with the inclusions. The property $\mathrm{tr}(ab) = \mathrm{tr}(ba)$ is a consequence of planar isotopy when all the strings added to $x$ to get $\widehat{x}$ go round $x$ in the same direction, as in Figure \ref{fig:tr(ab)=tr(ba)}.
Spherical isotopy reduces the general case to the one above.
Positive definiteness follows from the fact that the matrix units $e = \gamma \cdot \widetilde{\gamma'} \in P^{\mathcal{G}}_{\sigma\sigma^{\ast}}$ are mutually orthogonal elements of positive length: $\langle e, e \rangle = [3]^{-|\sigma|} \phi_{v_1} \phi_{v_2} > 0$, where $\gamma$, $\gamma'$ are paths of length $|\sigma|$ with edges on $\mathcal{G}$, $\mathcal{G}^{\mathrm{op}}$ as dictated by the sign string $\sigma$, which start at vertex $v_1$ and end at vertex $v_2$, and $\phi_v > 0$ for all $v$ since $\phi$ is a Perron-Frobenius eigenvector.
\hfill
$\Box$

\subsection{$P^{\mathcal{G}}$ as an $A_2$-$PTL$-module}

We will now determine all the irreducible weight-zero $A_2$-$ATL$-submodules of the $A_2$-graph planar algebra $P^{\mathcal{G}}$.
Let $\Delta_{\mathcal{G}} = G_{\rho}$ denote the adjacency matrix of the graph $\mathcal{G}$. If $\mathcal{G}$ is three-colourable then by permuting the labels $\Delta_{\mathcal{G}}$ may be written in the form
$$\Delta_{\mathcal{G}} = \left( { \begin{array}{ccc}
                 0 & \Delta_{01} & 0 \\
                 0 & 0 & \Delta_{12} \\
                 \Delta_{20} & 0 & 0
               \end{array} } \right),$$
where $\Delta_{01}$, $\Delta_{12}$ and $\Delta_{20}$ are matrices which give the number of edges between each 0,1,2-coloured vertex respectively of $\mathcal{G}$ to each 1,2,0-coloured vertex respectively. By a suitable ordering of the vertices the matrix $\Delta_{12}$ may be chosen to be symmetric. These matrices satisfy the conditions $\Delta_{01}^T \Delta_{01} = \Delta_{20} \Delta_{20}^T = \Delta_{12}^2$, $\Delta_{01} \Delta_{01}^T = \Delta_{20}^T \Delta_{20}$, which follow from the fact that $\Delta_{\mathcal{G}}$ is normal \cite[p.924]{evans/pugh:2009ii}. For non-three-colourable $\mathcal{G}$, we define $\Delta_{01} = \Delta_{12} = \Delta_{20} = \Delta_{\mathcal{G}}$.

Let $\beta_l^3$, $l \in \mathfrak{V}^{\mathcal{G}}_0$, be the eigenvalues of $\Delta_{01}\Delta_{12}\Delta_{20}$, and $v^{(l)}$ their corresponding eigenvectors. Then $(\Delta_{01}\Delta_{12}\Delta_{20})^T v^{(l)} = \overline{\beta}_l^3 v^{(l)}$ and $(\Delta_{01} \Delta_{01}^T)^3 v^{(l)} = \Delta_{01}\Delta_{12}\Delta_{20} (\Delta_{01}\Delta_{12}\Delta_{20})^T v^{(l)} = |\beta_l|^6 v^{(l)}$. Then if $\lambda_l$ are the eigenvalues of $\Delta_{01} \Delta_{01}^T$ with corresponding eigenvectors $v^{(l)}{}'$, $l \in \mathfrak{V}^{\mathcal{G}}_0$, we have $(\Delta_{01} \Delta_{01}^T)^3 v^{(l)}{}' = \lambda_l^3 v^{(l)}{}'$ so that $v^{(l)}{}' = v^{(l)}$ and $\lambda_l = |\beta_l|^2$. For the vertices $l \in \mathfrak{V}^{\mathcal{G}}_a$, where $a \in \{1,2\}$, we also have eigenvalues $\beta_l^3$ of $\Delta_{a,a+1}\Delta_{a+1,a+2}\Delta_{a+2,a}$.
However, comparing $\Delta_{20}\Delta_{01}\Delta_{12}v^{(l')} = \beta_{l'}^3 v^{(l')}$ where $l' \in \mathfrak{V}^{\mathcal{G}}_2$, and $\Delta_{20}\Delta_{01}\Delta_{12}\Delta_{20}v^{(l)} = \beta_l^3 \Delta_{20}v^{(l)}$ where $l \in \mathfrak{V}^{\mathcal{G}}_0$, we see that for each $\beta_{l_a} \neq 0$, where $l_a \in \mathfrak{V}^{\mathcal{G}}_a$, $a=1,2$, we have $\beta_{l_a}^3 = \beta_{l}^3$ for some $l \in \mathfrak{V}^{\mathcal{G}}_0$.

The dimension of $P^{\mathcal{G}}_{+^m -^n}$ is given by the trace of $\Delta^m (\Delta^T)^n$, which counts the number of pairs of paths on $\mathcal{G}$, $\mathcal{G}^{\mathrm{op}}$, and is given by the sum $\sum_l \beta_l^m \overline{\beta}_l^n$ of its eigenvalues, $l =1,2,\ldots, s$.
We can deduce all the irreducible weight-zero $A_2$-$ATL$-submodules of $P^{\mathcal{G}}$ in a similar way to \cite[Prop. 13]{reznikoff:2005}:

\begin{Prop} \label{Prop:RezProp13}
Let $\mathcal{G}$ be one of the finite $SU(3)$ $\mathcal{ADE}$ graphs, let $\zeta_l$ be the non-zero eigenvalues of $\Delta_{01}\Delta_{12}\Delta_{20}$, counting multiplicity, and let $\beta_l$ be any cubic root of $\zeta_l$, $l=1,2,\ldots,s'$ where $s = \textrm{min}\{s_0,s_1\}$. For all the three-colourable graphs except $\mathcal{E}_5^{(12)}$, we have $s_0 \geq s_1$, all the irreducible weight-zero $A_2$-$ATL$-submodules of the $A_2$-graph planar algebra $P^{\mathcal{G}}$ are $H^{\beta_l}$, $l=1,2,\ldots,s_1$, and $(s_0 - s_1)$ copies of $H^0$, and these can be assumed to be mutually orthogonal. For $\mathcal{E}_5^{(12)}$ we have $s_1 > s_0$, and all the irreducible weight-zero $A_2$-$ATL$-submodules of $P^{\mathcal{E}_5^{(12)}}$ are $H^{\beta_l}$, $l=1,2,\ldots,s_0$, and $2(s_1 - s_0)$ copies of $H^0$, which can again be assumed to be mutually orthogonal. If $\mathcal{G}$ is not three-colourable, all the irreducible weight-zero $A_2$-$ATL$-submodules of $P^{\mathcal{G}}$ are $H^{\beta_l}$, $l=1,2,\ldots,s_0$, where $s_0$ is the total number of vertices of $\mathcal{G}$.
\end{Prop}
\emph{Proof:}
Consider the case where $s_0 > s_1$ (the case for $\mathcal{E}_5^{(12)}$ where $s_1 > s_0$ is similar). For $l=1,2,\ldots,s'$, each $\beta_l$-eigenvector $v^{(l)} = (v^{(l)}_w)$, $w \in \mathfrak{V}^{\mathcal{G}}_0$ of $\Delta_{01} \Delta_{01}^T$ spans a one-dimensional subspace of $P^{\mathcal{G},\overline{0}}_{\varnothing}$ that is invariant under $A_2$-$ATL_{\varnothing}^{\overline{0}}$. To see this, first consider the element $\sigma_{01}\sigma_{12}\sigma_{20}$:
\begin{equation} \label{eqn:evector_invariance1}
\sigma_{01}\sigma_{12}\sigma_{20} v^{(l)} = \sigma_{01}\sigma_{12}\sigma_{20} \sum_{w \in \mathfrak{V}^{\mathcal{G}}_0} v^{(l)}_w = \sum_{w', w} (\Delta_{01} \Delta_{12} \Delta_{20})_{w',w} v^{(l)}_w = \sum_{w'} \beta_l^3 v^{(l)}_{w'} = \beta_l^3 v^{(l)},
\end{equation}
since $\beta_l^3$ is an eigenvalue for $\Delta_{01}\Delta_{12}\Delta_{20}$ for eigenvector $v^{(l)}$.
Similarly for $\sigma_{20}^{\ast}\sigma_{12}^{\ast}\sigma_{01}^{\ast}$. Next consider the general element $\sigma$ given by the composition of $2k$ elements $\sigma = \sigma_{01} \sigma_{12} \sigma_{20} \sigma_{01} \cdots \sigma_{\overline{k-1},\overline{k}} \sigma_{\overline{k-1},\overline{k}}^{\ast} \cdots \sigma_{12}^{\ast} \sigma_{01}^{\ast}$:
\begin{eqnarray}
\sigma v^{(l)} & = & \sum_{w', w} (\Delta_{01} \Delta_{12} \cdots \Delta_{\overline{k-1},\overline{k}} \Delta_{\overline{k-1},\overline{k}}^T \cdots \Delta_{01}^T)_{w',w} v^{(l)}_w \nonumber \\
& = & \sum_{w', w} ((\Delta_{01}\Delta_{01}^T)^k)_{w',w} v^{(l)}_w \;\; = \;\; \sum_{w'} |\beta_l|^{2k} v^{(l)}_{w'} \;\; = \;\; |\beta_l|^{2k} v^{(l)}. \label{eqn:evector_invariance2}
\end{eqnarray}
Any element of $A_2$-$ATL_{\varnothing}^{\overline{0}}$ is a linear combination of products of elements $\sigma_{j,j\pm1}$ such that the regions which meet the outer and inner boundaries have colour 0. Let $\sigma$ be such an element. Then the action of $\sigma$ on the $\beta_l$-eigenvector $v^{(l)}$ is given by $\sigma v^{(l)} = \sum_{w',w} M(w',w) v_w^{(l)}$, where $M$ is the product of matrices $\Delta_{\mathcal{G}}$, $\Delta_{\mathcal{G}}^T$ given by replacing every $\sigma_{j,j+1}$, $\sigma_{j',j'-1}$ in $\sigma$ by $\Delta_{\mathcal{G}}$, $\Delta_{\mathcal{G}}^T$ respectively. Then by (\ref{eqn:evector_invariance1}) and (\ref{eqn:evector_invariance2}), this gives some scalar multiple of $v^{(l)}$. Similarly $\sigma_{a,0} v^{(l)}$ spans a one-dimensional subspace of $P^{\mathcal{G},\overline{a}}_{\varnothing}$ that is invariant under $A_2$-$ATL_{\varnothing}^{\overline{a}}$, $a=1,2$.
Then for each $l=1,2,\ldots,s'$, the $\beta_l$-eigenvector $v^{(l)}$ generates the submodule $H^{\beta_l}$ by Proposition \ref{Thm:SU(3)5.12}. The inner product on $H^{\beta_l}$ coincides with the inner product on $P^{\mathcal{G}}$. To see this we only need to check its restriction to the zero-weight part because of (\ref{eqn:inner-product_on_Hilbert-module}). For any element $A \in A_2 \mathrm{-} ATL_{\varnothing}^{\overline{0}}$, $\langle Av, v \rangle_{H^{\beta_l}} = c \langle v, v \rangle_{H^{\beta_l}}$ whilst $\langle A v^{(l)}, v^{(l)} \rangle_{P^{\mathcal{G}}} = d \langle v^{(l)}, v^{(l)} \rangle_{P^{\mathcal{G}}}$. The element $A$ is necessarily a combination of non-contractible circles, which gives the same contribution in $P^{\mathcal{G}}$ as in $H^{\beta_l}$ by (\ref{eqn:evector_invariance1}), (\ref{eqn:evector_invariance2}). So $c = d$. This shows that the inner product on the $H^{\beta_l}$ is positive definite, since the inner product on $P^{\mathcal{G}}$ is.

Similarly, a 0-eigenvector generates the submodule $H^0$, where for $s_0 > s_1$, $\mathrm{dim}(H^{0,\overline{0}}_{\varnothing,\overline{0}}) = 1$ and $\mathrm{dim}(H^{0,\overline{1}}_{\varnothing,\overline{1}}) = \mathrm{dim}(H^{0,\overline{2}}_{\varnothing,\overline{2}}) = 0$, whilst for $\mathcal{E}_5^{(12)}$ we have $\mathrm{dim}(H^{0,\overline{1}}_{\varnothing,\overline{1}}) = \mathrm{dim}(H^{0,\overline{2}}_{\varnothing,\overline{2}}) = 1$ and $\mathrm{dim}(H^{0,\overline{0}}_{\varnothing,\overline{0}}) = 0$. As in \cite{reznikoff:2005}, in order to make the resulting submodules orthogonal we take an orthogonal set of eigenvectors.
\hfill
$\Box$

For an $\mathcal{ADE}$ graph $\mathcal{G}$ with Coxeter number $n$, let $\beta_{(l_1,l_2)}$ be the eigenvalue given by
\begin{equation} \label{eqn:evalues-A(n)}
\beta_{(l_1,l_2)} = \exp \left( \frac{2i \pi}{3n} (l_1 + 2l_2 +3) \right) + \exp \left( - \frac{2i \pi}{3n} (2l_1 + l_2 +3) \right) + \exp \left( \frac{2i \pi}{3n} (l_1 - l_2) \right)
\end{equation}
for exponent $(l_1,l_2)$ of $\mathcal{G}$ (see, e.g. \cite{evans/pugh:2010i}).
Then for the graphs $\mathcal{A}^{(n)}$,
\begin{equation}
P^{\mathcal{A}^{(n)}} \supset \bigoplus_{(l_1,l_2)} H^{\beta_{(l_1,l_2)}},
\end{equation}
for $n \not \equiv 0 \textrm{ mod } 3$, whilst for $n=3k$, $k \geq 2$,
\begin{equation}
P^{\mathcal{A}^{(3k)}} \supset \bigoplus_{(l_1,l_2)} H^{\beta_{(l_1,l_2)}} \oplus H^{0,\overline{0}},
\end{equation}
where in both cases the summation is over all $(l_1,l_2) \in \{ (m_1,m_2) | \; 3m_2 \leq n-3, 3m_1 + 3m_2 < 2n-6 \}$, i.e. each $\beta_{(l_1,l_2)}$ is a cubic root of an eigenvalue of $\Delta_{01} \Delta_{12} \Delta_{20}$.
We believe that we in fact have equality here, so that $P^{\mathcal{A}^{(n)}} = \bigoplus_{(l_1,l_2)} H^{\beta_{(l_1,l_2)}}$. In the $A_1$-case this was achieved by a dimension count of the left and right hand sides \cite[Thm. 15]{reznikoff:2005}. However, we have not yet been able to determine a similar result in our $A_2$-setting.

For the other $\mathcal{ADE}$ graphs, Prop. \ref{Prop:RezProp13} gives the following results for the zero-weight part of $P^{\mathcal{G}}$. For the $\mathcal{D}$ graphs, we have
\begin{equation}
P^{\mathcal{D}^{(3k)}} \supset \bigoplus_{(l_1,l_2)} H^{\beta_{(l_1,l_2)}} \oplus 3H^{0,\overline{0}},
\end{equation}
for $k \geq 2$, where the summation is over all $(l_1,l_2) \in \{ (m_1,m_2) | m_2 \leq k-1, m_1 + m_2 < 2k-2, m_1 - m_2 \equiv 0 \textrm{ mod } 3 \}$, whilst for $n \not \equiv 0 \textrm{ mod } 3$,
\begin{equation}
P^{\mathcal{D}^{(n)}} \supset \bigoplus_{(l_1,l_2)} H^{\beta_{(l_1,l_2)}},
\end{equation}
where the summation is over all $(l_1,l_2) \in \{ (m_1,m_2) | 3m_2 \leq n-3, 3m_1 + 3m_2 < 2n-6 \}$.
The path algebras for $\mathcal{A}^{(n)\ast}$ and $\mathcal{D}^{(n)\ast}$ are identified under the map which send the vertices $i_l$, $j_l$ and $k_l$ of $\mathcal{D}^{(n)\ast}$ with the vertex $l$ of $\mathcal{A}^{(n)\ast}$, $l=1,2,\ldots, \lfloor l/2 \rfloor$. We have
\begin{equation}
P^{\mathcal{A}^{(n)\ast}} = P^{\mathcal{D}^{(n)\ast}} \supset \bigoplus_{(l_1,l_2)} H^{\beta_{(l_1,l_2)}},
\end{equation}
where the summation is over all $(l_1,l_2) \in \{ (m,m) | m=0,1,\ldots, \lfloor (n-3)/2 \rfloor \}$. Similarly, the path algebras for $\mathcal{E}^{(8)}$ and $\mathcal{E}^{(8)\ast}$ are identified, and
\begin{equation}
P^{\mathcal{E}^{(8)}} = P^{\mathcal{E}^{(8)\ast}} \supset H^{\beta_{(0,0)}} \oplus H^{\beta_{(3,0)}} \oplus H^{\beta_{(0,3)}} \oplus H^{\beta_{(2,2)}}.
\end{equation}
For the graphs $\mathcal{E}_i^{(12)}$, $i=1,2,3$, we have
\begin{equation}
P^{\mathcal{E}_i^{(12)}} \supset H^{\beta_{(0,0)}} \oplus 2 H^{\beta_{(2,2)}} \oplus H^{\beta_{(4,4)}}.
\end{equation}
For the remaining exceptional graphs we have
\begin{eqnarray}
P^{\mathcal{E}_4^{(12)}} & \supset & H^{\beta_{(0,0)}} \oplus H^{\beta_{(2,2)}} \oplus H^{\beta_{(4,4)}} \oplus 2H^{0,\overline{0}}, \\
P^{\mathcal{E}_5^{(12)}} & \supset & H^{\beta_{(0,0)}} \oplus H^{\beta_{(3,0)}} \oplus H^{\beta_{(0,3)}} \oplus H^{\beta_{(2,2)}} \oplus H^{\beta_{(4,4)}} \oplus H^{0,\overline{1}} \oplus H^{0,\overline{2}}, \\
P^{\mathcal{E}^{(24)}} & \supset & H^{\beta_{(0,0)}} \oplus H^{\beta_{(6,0)}} \oplus H^{\beta_{(0,6)}} \oplus H^{\beta_{(4,4)}} \oplus H^{\beta_{(7,4)}} \oplus H^{\beta_{(4,7)}} \oplus H^{\beta_{(6,6)}} \oplus H^{\beta_{(10,10)}}. \qquad
\end{eqnarray}

The $A_2$-planar algebra $P \cong PTL$ associated to the subfactor double complex for $\mathcal{A}^{(n)}$ \cite{evans/pugh:2009iii}, clearly has decomposition $P = H^{\alpha}$ as an $A_2$-$ATL$-module, since $PTL$ is equal to the $A_2$-$ATL$-module $H^{\alpha}$ (see Section \ref{sect:Irr_PTL-modules}). Since every $A_2$-planar algebra contains $PTL$, the $A_2$-$C^{\ast}$-planar algebra associated to the subfactor double complex for any $\mathcal{ADE}$ graph with a flat connection will contain the zero-weight module $H^{\alpha}$. The graphs $\mathcal{A}^{(n)}$, $\mathcal{D}^{(n)}$, $n < \infty$, have a flat connection \cite{evans/kawahigashi:1994}, and it is expected that the only other graphs with a flat connection are $\mathcal{E}^{(8)}$, $\mathcal{E}_1^{(12)}$ and $\mathcal{E}^{(24)}$ \cite[pp.359-360]{evans/pugh:2009iii}.

\subsection{Irreducible modules with non-zero weight}

In this section we present some irreducible $A_2$-$ATL$-modules with non-zero weight.

\begin{figure}[tb]
\begin{minipage}[t]{9cm}
\begin{center}
  \includegraphics[width=30mm]{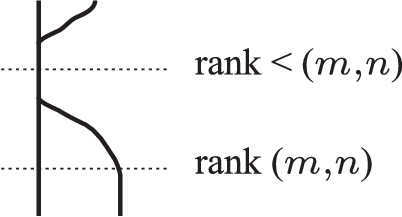}
 \caption{Two paths sharing a common string} \label{fig-through_strings_common}
\end{center}
\end{minipage}
\hfill
\begin{minipage}[t]{6cm}
\begin{center}
  \includegraphics[width=25mm]{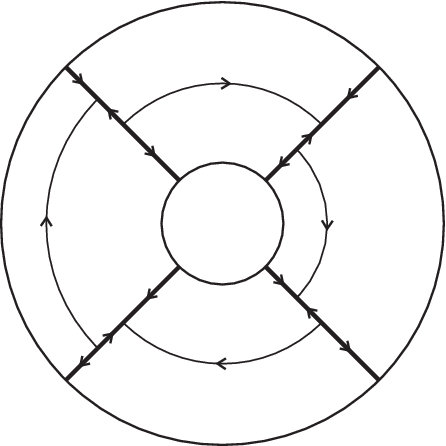}
 \caption{A tangle in $A^{(2,2)}$} \label{fig-tangle-A(2,2)}
\end{center}
\end{minipage}
\end{figure}

Let $A^{(m,n)} := A_2\textrm{-}ATL_{+^m -^n}/A_2\textrm{-}ATL_{+^m -^n}^{(m,n)}$ and $A^{(m,n)}_0 := A_2\textrm{-}ATL_{+^m -^n}^{\overline{0}}/A_2\textrm{-}ATL_{+^m -^n}^{\overline{0},(m,n)}$, where $A_2\textrm{-}ATL_{+^m -^n}^{\overline{0},(m,n)}$ is the two-sided ideal in $A_2\textrm{-}ATL_{+^m -^n}^{\overline{0}}$ given by the linear span of all labeled $A_2$-annular $+^m -^n$-tangles in $A_2\textrm{-}ATL_{+^m -^n}^{\overline{0}}$ with rank $< (m,n)$.
A \textbf{through string} in a tangle in $A^{(m,n)}$ will be a choice of a path along strings from a vertex on the outer disc of the tangle to a vertex on its inner disc, such that this choice allows each vertex on the outer and inner disc to be the endpoint of disjoint through strings, that is, no two through strings share a common string. It is always possible to make such a choice of $m+n$ through strings in any tangle $T$ in $A^{(m,n)}$. To see this, note that the rank at the inner and outer disc of $T$ is $(m,n)$, thus if two paths from the outer disc to the inner disc share a common string, the rank of $T$ would necessarily be less than $(m,n)$, as illustrated in Figure \ref{fig-through_strings_common}, and thus $T$ is in the ideal $A_2\textrm{-}ATL_{+^m -^n}^{(m,n)}$.

Thus we will draw any basis tangle in $A^{(m,n)}$ with $m+n$ through strings (denoted by thick strings). These through strings will be connected by other strings, which must be single strings which do not have any trivalent vertices along them otherwise digons or embedded squares must necessarily be created- we will call these connecting strings \textbf{rungs}. An example of a tangle in $A^{(2,2)}$ is illustrated in Figure \ref{fig-tangle-A(2,2)}.

\subsubsection{The case $m,n \neq 0$}

For $m,n \neq 0$, we will show that there is a unique choice of through strings in any tangle in $A^{(m,n)}$. First we need some more definitions and a lemma. We will call a region in a tangle a \textbf{through space} if its boundary contains a segment of both the inner and outer disc, e.g. the regions $R,S$ are through spaces in the following tangle:
\begin{center}
\includegraphics[width=25mm]{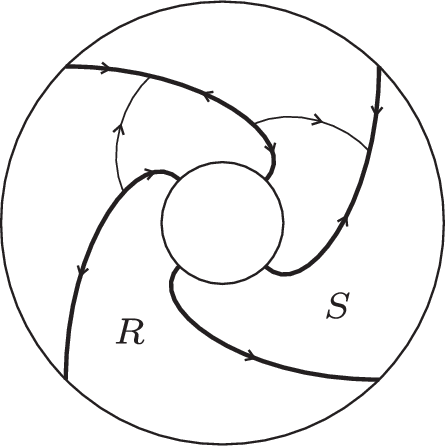}
\end{center}

\begin{figure}[htb]
\begin{center}
  \includegraphics[width=20mm]{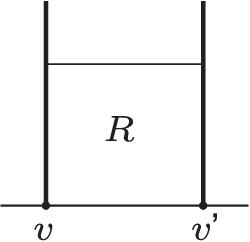}
 \caption{A bottom rung} \label{fig-bottom_rung}
\end{center}
\end{figure}

If there is a region $R$ whose boundary contains a segment of the outer or inner discs, and the rest of its boundary consists of three strings joined together at two trivalent vertices (see Figure \ref{fig-bottom_rung}), we will call the string between the two trivalent vertices a \textbf{bottom rung}.

\begin{figure}[htb]
\begin{center}
 \begin{minipage}[b]{6cm}
  \begin{center}
  \includegraphics[width=25mm]{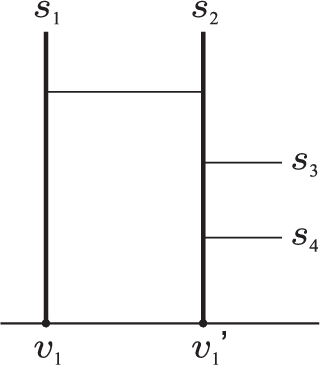}
  \end{center}
 \end{minipage}
 \begin{minipage}[b]{6cm}
  \begin{center}
  \includegraphics[width=30mm]{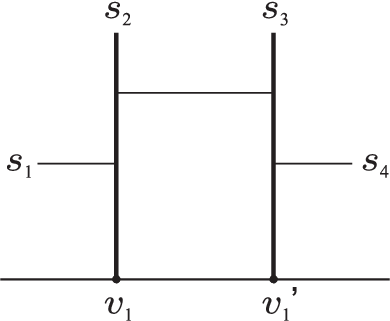}
  \end{center}
 \end{minipage} \\
 \caption{Choices of through strings for a region $R$ bounded by five strings} \label{fig-(in)feasible_choice}
\end{center}
\end{figure}

\begin{Lemma} \label{Lemma:through_space/rung}
Any basis tangle in $A^{(m,n)}$, where $m,n \neq 0$, contains either a through space or a bottom rung.
\end{Lemma}
\emph{Proof:}
Let $T$ be any basis tangle in $A^{(m,n)}$. There are two pair of adjacent vertices $(v_j,v_j')$ on the outer disc which have opposite orientations, i.e. $v_j = \pm$, $v_j' = \mp$, where $j=1,2$.
Suppose $T$ does not contain any through spaces, and consider the region $R$ which meets the outer disc at the segment between vertices $v_1$ and $v_1'$. If this region is only bounded by three strings, as in Figure \ref{fig-bottom_rung}, then $T$ contains a bottom rung. Now suppose the region $R$ is bounded by five strings:
\begin{center}
\includegraphics[width=40mm]{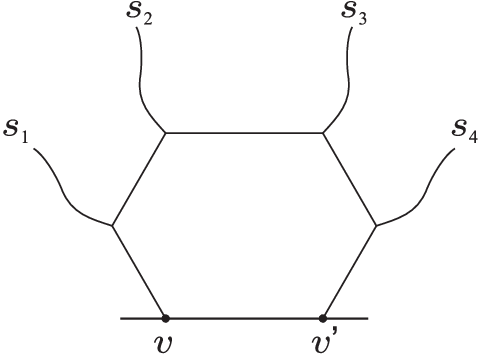}
\end{center}
The through strings cannot be both $s_1$ and $s_2$ (and similarly for both $s_3$ and $s_4$), as in the left-hand side of Figure \ref{fig-(in)feasible_choice}, since then both $s_3$ and $s_4$ must be connected to the adjacent through string, which creates an embedded, contradicting the assumption that $T$ is a basis tangle. A similar argument can also be used to show that this region $R$ cannot be bounded by more than five strings in a basis tangle. Thus the through strings must be $s_2$ and $s_3$, as illustrated on the right-hand side of Figure \ref{fig-(in)feasible_choice}.

Recall that $v_1$, $v_1'$ have opposite orientations. Due to the orientations of the adjacent boundary points, we must have the following situation:
\begin{center}
\includegraphics[width=50mm]{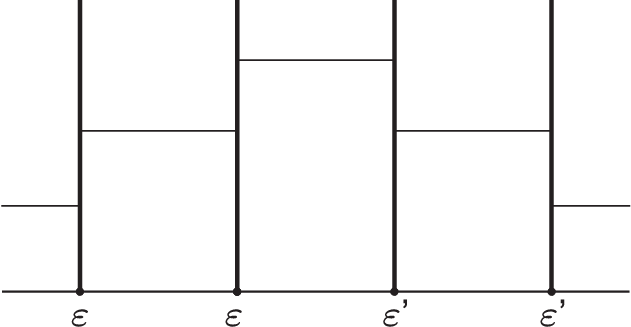}
\end{center}
where the orientations of the boundary points are $\varepsilon = v_1$, $\varepsilon' = v_1'$, and thus at $(v_2,v_2')$ we must have a bottom rung:
\begin{center}
\includegraphics[width=35mm]{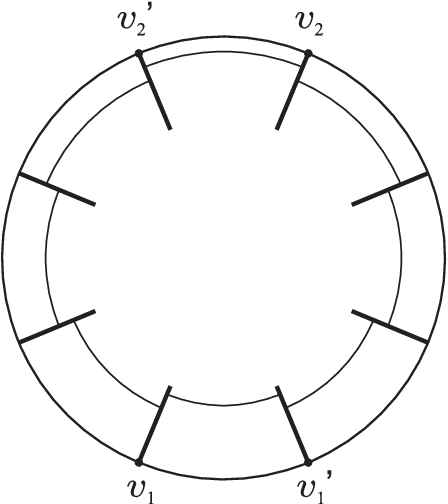}
\end{center}
\hfill
$\Box$

A similar argument gives the following:
\begin{Cor} \label{Cor:through_space/rung}
Let $\sigma$ be any permutation of $+^m -^n$, for any $m,n \neq 0$. Then any basis $((\sigma,\overline{a}),(+^m -^n,\overline{a'}))$-tangle of rank $(m,n)$ contains either a through string of a bottom rung.
\end{Cor}

Then we have the following result:
\begin{Lemma}
When $m,n \neq 0$, there is only one possible choice of $m+n$ disjoint through strings in any basis tangle in $A^{(m,n)}$.
\end{Lemma}
\emph{Proof:}
Let $T$ be any basis tangle in $A^{(m,n)}$. Suppose first that $T$ contains a through space. Since the rungs in $T$ must connect two distinct through strings, the strings which form the boundaries of the through space must be a pair of adjacent through strings, and by the same argument the other through strings are uniquely determined.

Now suppose that the basis tangle $T$ does not contain a through space. Then by Lemma \ref{Lemma:through_space/rung}, it must contain a bottom rung. Here there is only one possibility for the allocation of disjoint through strings:
\begin{center}
\includegraphics[width=20mm]{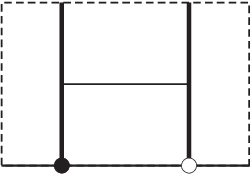}
\end{center}
Here the vertex $\bullet$ has the opposite orientation to the vertex $\circ$. We now consider the tangle $T'$ obtained by replacing the part of $T$ contained within the dashed rectangle above by
\begin{center}
\includegraphics[width=20mm]{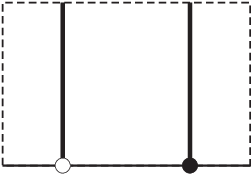}
\end{center}
and leaving the rest of the tangle unchanged. Note that $T'$ is no longer a $+^m -^n$-tangle as the orientations of the two vertices have been reversed. However, by Corollary \ref{Cor:through_space/rung}, $T'$ must either contain a through space or a bottom rung, and so we proceed as we did for the original tangle $T$. This procedure reduces the number of rungs in the tangle by one, and thus iterating this procedure will lead to a tangle which contains a through space. Thus we find the unique way of choosing the through strings in $T$.
\hfill
$\Box$

\begin{figure}[htb]
\begin{center}
  \includegraphics[width=120mm]{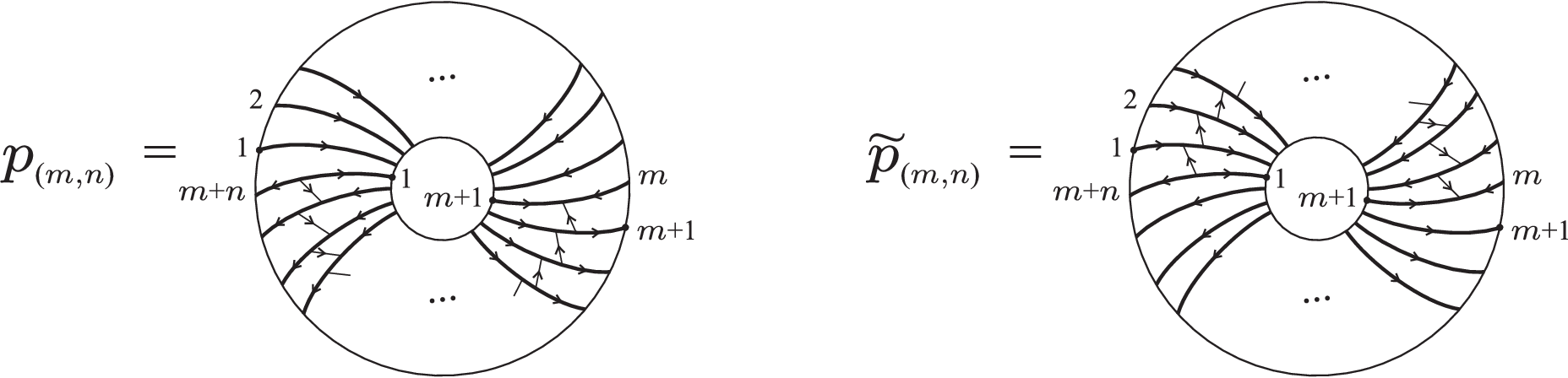}
 \caption{$+^m -^n$-tangles $p_{(m,n)}$ and $\widetilde{p}_{(m,n)}$} \label{fig-p(m,n)}
\end{center}
\end{figure}

We will now show that any basis tangle in $A^{(m,n)}$ can be written as $p_{(m,n)}^i \widetilde{p}_{(m,n)}^j$, $i,j \in \mathbb{Z}$, where the $+^m -^n$-tangles $p_{(m,n)}$, $\widetilde{p}_{(m,n)}$ are illustrated in Figure \ref{fig-p(m,n)}. Note that $p_{(m,n)}^{\ast} = p_{(m,n)}^{-1}$, $\widetilde{p}_{(m,n)}^{\ast} = \widetilde{p}_{(m,n)}^{-1}$. We will usually simply write $p$, $\widetilde{p}$ for $p_{(m,n)}$, $\widetilde{p}_{(m,n)}$ respectively, where $(m,n)$ is obvious. It is easy to see by drawing pictures, and using relations K1-K3 if necessary, that $p$, $\widetilde{p}$, $p^{\ast}$, $\widetilde{p}^{\ast}$ all commute in $A^{(m,n)}$. Note that if $T$ is a $+^m -^n$-tangle with outer disc labeled by $\ast_a$, then $pT$, $\widetilde{p}T$ are $+^m -^n$-tangles with outer disc labeled by $\ast_{a+1}$, $\ast_{a-1}$ respectively, i.e. for some $A_2\textrm{-}ATL$-module $V$, $p$, $\widetilde{p}$ is a map from $V_{+^m -^n}^{\overline{a}}$ to $V_{+^m -^n}^{\overline{a+1}}$, $V_{+^m -^n}^{\overline{a-1}}$ respectively.
The tangles $p$, $\widetilde{p}$ can be regarded as some sort of ``rotation by one'', since the through strings connect vertex $i$ on the outer disc to vertex $i+1$ on the inner disc. We will say that a basis tangle $T$ in $A^{(m,n)}$ has \textbf{rotation number} $b \in \mathbb{N}\cup\{0\}$ if the through strings connect vertex $i$ on the outer disc to vertex $i+b$ on the inner disc.

By the same argument as was used at the end of the proof of Lemma \ref{Lemma:through_space/rung}, any basis tangle $T$ in $A^{(m,n)}$ which does not contain a through space may be written as $T=XT'$, where $X$ is a tangle with exactly one rung between each pair of adjacent through strings, as illustrated in Figure \ref{fig-tangle-A(2,2)}. Such a tangle is given by $p\widetilde{p}^{\ast}$ or its inverse. Repeating this argument for $T'$, we see that any basis tangle $T$ in $A^{(m,n)}$ which does not contain a through space may be written as $T=X^i T'$, where $X = p\widetilde{p}^{\ast}$, $i \in \mathbb{Z}$, and $T'$ is a basis tangle in $A^{(m,n)}$ which contains a through space. Thus to shows that any basis tangle in $A^{(m,n)}$ can be written as $p^i \widetilde{p}^j$ we need only show that any basis tangle $T'$ which contains a through space can be written in this form. We will show this by induction on the rotation number of such $T'$.

We start with rotation number 0. If there is a rung anywhere in $T'$, then the orientations of the vertices on the inner and outer discs dictate that we must have a factor $p\widetilde{p}^{\ast}$, which contradicts the fact that $T'$ contains a through space. Hence the only basis tangle $T'$ with rotation number 0 which contains a through space is the identity tangle $\mathbf{1}_{+^m -^n}$.

Suppose now that for any basis tangle which contains a through space and has rotation number $k$ can be written as $p^i \widetilde{p}^j$ for some $i,j \in \mathbb{Z}$, and let $T'$ be any basis tangle which contains a through space and has rotation number $k+1$. Then $T' p^{\ast}$ has rotation number $k$, and thus $T' p^{\ast} = p^i \widetilde{p}^j$ for some $i,j \in \mathbb{Z}$. Then $T' = p^{i+1} \widetilde{p}^j$ as required. Thus we have proved the following:

\begin{Lemma}
Any basis tangle in $A^{(m,n)}$ can be written as $p^i \widetilde{p}^j$ for some $i,j \in \mathbb{Z}$.
\end{Lemma}

\begin{figure}[htb]
\begin{center}
  \includegraphics[width=120mm]{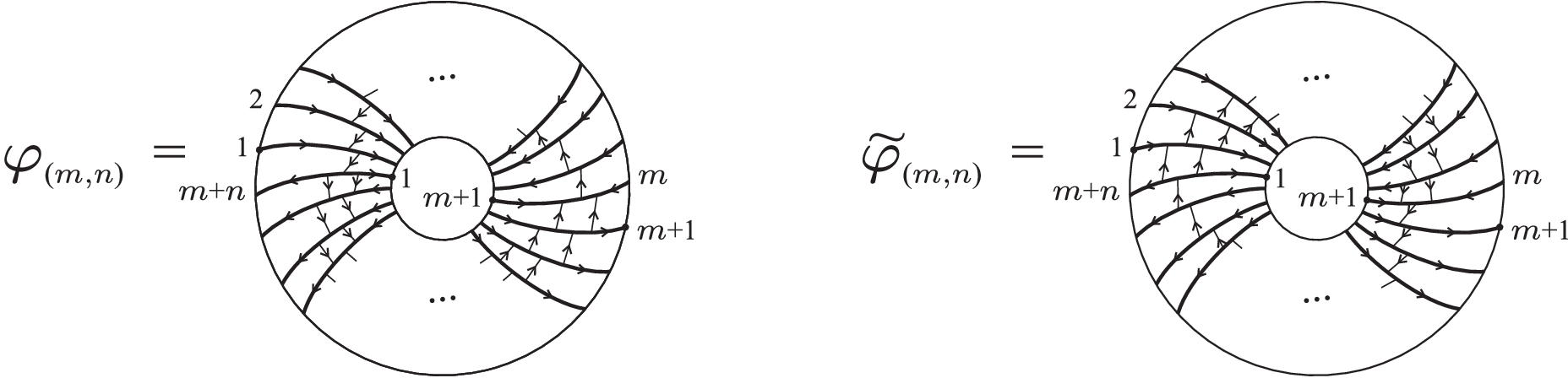}
 \caption{$+^m -^n$-tangles $\varphi_{(m,n)}$ and $\widetilde{\varphi}_{(m,n)}$} \label{fig-varphi(m,n)}
\end{center}
\end{figure}

Now $A^{(m,n)}_0 = \{ p^i \widetilde{p}^j | \; i-j \equiv 0 \textrm{ mod } 3 \}$, thus $A^{(m,n)}_0$ is generated by $p^3$, $\widetilde{p}^3$, $p^2 \widetilde{p}^{\ast}$ and $p^{\ast} \widetilde{p}^2$. We let $\varphi_{(m,n)} := p_{(m,n)}^2 \widetilde{p}_{(m,n)}^{\ast}$ and $\widetilde{\varphi}_{(m,n)} := p_{(m,n)}^{\ast} \widetilde{p}_{(m,n)}^2$. These tangles are illustrated in Figure \ref{fig-varphi(m,n)}. Then since $\varphi^2 \widetilde{\varphi} = p^4 (\widetilde{p}^{\ast})^2 p^{\ast} \widetilde{p}^2 = p^3$, and $\varphi \widetilde{\varphi}^2 = \widetilde{p}^3$, we see that $A^{(m,n)}_0$ is generated by $\varphi_{(m,n)}$ and $\widetilde{\varphi}_{(m,n)}$. Note that we have $\varphi^{\ast} = \varphi^{-1}$ and $\widetilde{\varphi}^{\ast} = \widetilde{\varphi}^{-1}$.

Let $\rho_{(m,n)}$ be the $+^m -^n$-tangle given by the image of $\varphi_{(m,n)} \widetilde{\varphi}_{(m,n)}$ in $A^{(m,n)}_0$, illustrated in Figure \ref{fig-rho(m,n)}. The tangle $\rho_{(m,n)}$ is a ``rotation by two''. Indeed, it is easy to see by drawing pictures and using relation K3 that $\rho_{(m,n)}$ is a rotation of order $k$ in $A^{(m,n)}_0$, i.e. $(\rho_{(m,n)})^k = \mathbf{1}_{+^m -^n}$ in $A^{(m,n)}_0$, where $k = m+n$ if $m+n$ is odd and $k=(m+n)/2$ if $m+n$ is even.

\begin{figure}[htb]
\begin{center}
  \includegraphics[width=60mm]{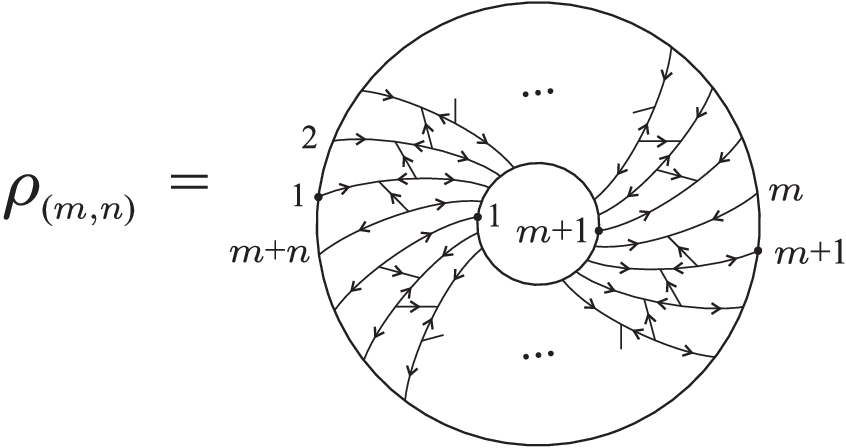}
 \caption{$+^m -^n$-tangle $\rho_{(m,n)}$} \label{fig-rho(m,n)}
\end{center}
\end{figure}

Note that the algebras $A^{(m,n)}$ are infinite dimensional since $p^i \widetilde{p}^j$ are all distinct tangles in $A^{(m,n)}$, for $i,j \in \mathbb{Z}$. However, we can obtain a finite-dimensional $A_2\textrm{-}ATL$-module $V^{(m,n),\gamma}$ by imposing the condition that $(p\widetilde{p})^3$ acts by multiplication by a scalar $\gamma \in \mathbb{C}$ in the lowest weight module $V^{(m,n),\gamma}_{+^m -^n}$, i.e. $(p\widetilde{p})^3 = \gamma \mathbf{1}_{+^m -^n}$ in $V^{(m,n),\gamma}$. Note that the $A_2\textrm{-}ATL$-module $V^{(m,n),\gamma}$ may not be irreducible, but may rather decompose into irreducible modules.
Now $\varphi^{m+n} = (p\widetilde{p})^{m+2n}$ and $\widetilde{\varphi}^{m+n} = (p\widetilde{p})^{2m+n}$. We have $\varphi^{m+n} = \gamma^{(m+2n)/3} \mathbf{1}$ and $\widetilde{\varphi}^{m+n} = \gamma^{(2m+n)/3} \mathbf{1}$ (note that $m+2n \equiv 2m+n \equiv 0 \textrm{ mod } 3$ since $m \equiv n \textrm{ mod } 3$). Then we have $\varphi^{\ast} = \gamma^{-(m+2n)/3} \varphi^{m+n-1}$ and $\widetilde{\varphi}^{\ast} = \gamma^{-(2m+n)/3} \widetilde{\varphi}^{m+n-1}$. Substituting for $\varphi^{\ast}$ in $\varphi^{\ast} \widetilde{\varphi} = ((p\widetilde{p})^3)^{\ast} = \overline{\gamma} \mathbf{1}$ we obtain $\varphi^{m+n-1} \widetilde{\varphi} = \overline{\gamma} \gamma^{(m+2n)/3} \mathbf{1}$, which gives the relation $\widetilde{\varphi}^{\ast} = \overline{\gamma}^{-1} \gamma^{-(m+2n)/3} \varphi^{m+n-1}$.
Then we have $\widetilde{\varphi} = \gamma^{-1} \overline{\gamma}^{-(m+2n)/3} (\varphi^{\ast})^{m+n-1}$, and since $\widetilde{\varphi}^{\ast} = \widetilde{\varphi}^{-1}$, we have $\widetilde{\varphi} \widetilde{\varphi}^{\ast} = |\gamma|^{-2} |\gamma|^{-2(m+2n)/3} \mathbf{1} = \mathbf{1}$, and thus we must take $\gamma \in \mathbb{T} = \{ \gamma \in \mathbb{C} : \; |\gamma| = 1\}$.
Then as an $A_2$-$ATL_{+^m -^n}^{\overline{0}}$-module
$$V^{(m,n),\gamma}_{+^m -^n} = \mathrm{span}(\varphi_{(m,n)}^l | \; l = 0,1,\ldots,m+n-1),$$
where $\varphi_{(m,n)}^{m+n} = \gamma^{(m+2n)/3} \mathbf{1}_{+^m -^n}$.
For all $m,n \neq 0$, we can choose a faithful trace $\mathrm{tr}'$ on $A^{(m,n)}$, which we extend to a trace $\mathrm{tr}$ on $A_2$-$ATL_{+^m -^n}$ by $\mathrm{tr} = \mathrm{tr}' \circ \pi$, where $\pi$ is the quotient map $\pi:A_2\textrm{-}ATL_{+^m -^n} \rightarrow A^{(m,n)}$. We can define an inner product on $A_2$-$ATL(\sigma:+^m -^n)$ by $\langle S,T \rangle = \mathrm{tr}(T^{\ast}S)$ for any $S,T \in A_2\textrm{-}ATL(\sigma:+^m -^n)$.

If $\gamma^{(m+2n)k/3} \neq 1$ for any $k \in \mathbb{N}$, then $\mathrm{dim}(V^{(m,n),\gamma}_{+^m -^n}) = m+n$.
The $A_2$-$ATL$-module $V^{(m,n),\gamma}$ is irreducible by Lemma \ref{Lemma:SU(3)Lemma3.4} since $V^{(m,n),\gamma} = A_2\mathrm{-}ATL(V^{(m,n),\gamma}_{+^m -^n})$.
If the inner product is positive semi-definite, we define the Hilbert $A_2$-$ATL$-module $H^{(m,n),\gamma}$, $\gamma^{(m+2n)k/3} \neq 1$ for any $k \in \mathbb{N}$, to be the quotient of $V^{(m,n),\gamma}$ by the zero-length vectors with respect to this inner product; otherwise $H^{(m,n),\gamma}$ does not exist.

Now suppose $\gamma^{(m+2n)/3} = 1$. Then we see that $\varphi_{(m,n)}$ acts on $V^{(m,n),\gamma}_{+^m -^n}$ as $\mathbb{Z}_{m+n}$, by permuting the $m+n$ basis elements $\varphi_{(m,n)}^l$, and so the $A_2$-$ATL_{+^m -^n}^{\overline{0}}$-module $V^{(m,n),\gamma}_{+^m -^n}$ decomposes as a direct sum over the $(m+n)^{\textrm{th}}$ roots of unity $\omega$ of $A_2$-$ATL_{+^m -^n}^{\overline{0}}$-modules $V^{(m,n),\gamma, \omega}_{+^m -^n}$, where $V^{(m,n),\gamma, \omega}_{+^m -^n}$ is the $\omega$-eigenspace for the action of $\mathbb{Z}_{m+n}$ with eigenvalue $\omega$.
Since $\varphi_{(m,n)}^{\ast} \varphi_{(m,n)} = \mathbf{1}_{+^m -^n}$, the decomposition into $V^{(m,n),\gamma, \omega}_{+^m -^n}$ is orthogonal. If we let $\psi_{(m,n)}^{\gamma,\omega}$ be the vector in $V^{(m,n),\gamma,\omega}_{+^m -^n}$ which is proportional to $\sum_{j=0}^{m+n-1} (\omega \gamma)^{-j} \varphi_{(m,n)}^j$ such that $\langle \psi_{(m,n)}^{\gamma,\omega}, \psi_{(m,n)}^{\gamma,\omega} \rangle = 1$, then $\varphi_{(m,n)} \psi_{(m,n)}^{\gamma,\omega} = \omega \gamma \psi_{(m,n)}^{\gamma,\omega}$. We see that $\mathrm{dim}(V^{(m,n),\gamma,\omega}_{+^m -^n}) = 1$, and $V^{(m,n),\gamma,\omega}_{+^m -^n}$ is the span of $\psi_{(m,n)}^{\gamma,\omega}$.
The $A_2$-$ATL$-module $V^{(m,n),\gamma,\omega}$ is irreducible by Lemma \ref{Lemma:SU(3)Lemma3.4} since $V^{(m,n),\gamma,\omega} = A_2\mathrm{-}ATL(\psi_{(m,n)}^{\gamma,\omega})$.

If the inner product is positive semi-definite, we define the Hilbert $A_2$-$ATL$-module $H^{(m,n),\gamma,\omega}$, $\gamma^{(m+2n)/3} = 1$, to be the quotient of $V^{(m,n),\gamma,\omega}$ by the zero-length vectors with respect to this inner product; otherwise $H^{(m,n),\gamma,\omega}$ does not exist.

\subsubsection{The case $n = 0$}

\begin{figure}[htb]
\begin{center}
  \includegraphics[width=140mm]{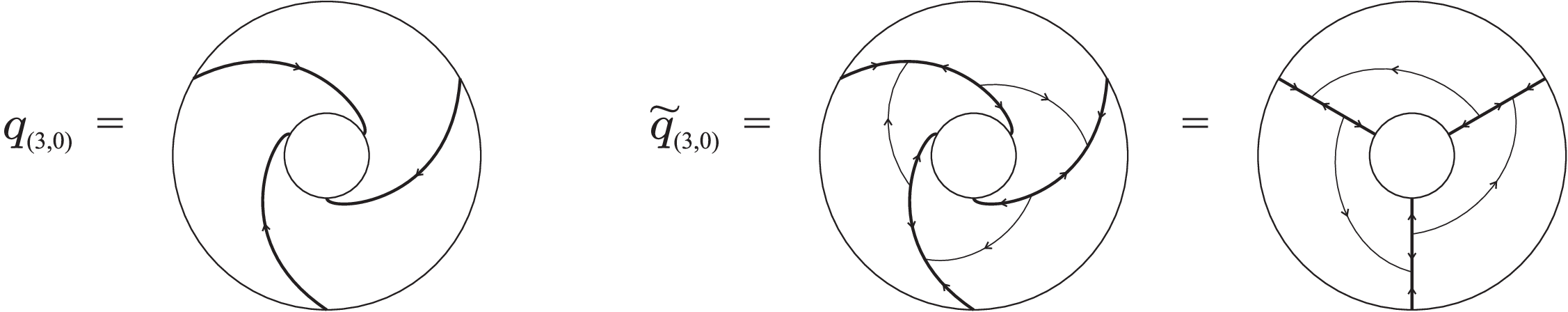}
 \caption{$+^3$-tangles $q_{(3,0)}$ and $\widetilde{q}_{(3,0)}$} \label{fig-q(3,0)}
\end{center}
\end{figure}

We will now consider the case when $n=0$. The case when $m=0$ is similar. When $n=0$, a basis tangle in $A^{(m,0)}$ will not always have a unique rotation number, as the choice of through strings is not always unique. Consider the case $m=3$. The algebra $A^{(3,0)}$ is generated by the tangles $q_{(3,0)}$ and $\widetilde{q}_{(3,0)}$, which are illustrated in Figure \ref{fig-q(3,0)}. Notice that for $\widetilde{q}_{(3,0)}$, the first choice of through strings gives rotation number 1, whilst the second choice gives rotation number 0. We will usually write $q$, $\widetilde{q}$ for $q_{(3,0)}$, $\widetilde{q}_{(3,0)}$. Again, it is easy to see by drawing pictures, and using relations K1-K3 if necessary, that $q$, $\widetilde{q}$, $q^{\ast}$, $\widetilde{q}^{\ast}$ all commute in $A^{(3,0)}$, and that $\widetilde{q}^{\ast} = q^{\ast} \widetilde{q}$.
Note that $q^{\ast} = q^{-1}$, but that $\widetilde{q}$, $\widetilde{q}^{\ast}$ do not have an inverse in $A^{(3,0)}$. Due to the orientations of the vertices on the inner and outer discs, there must be an equal number of rungs between any pair of adjacent through strings.
Both $q$ and $\widetilde{q}$ can be drawn as tangles with rotation number 1. The only other basis tangle with rotation number 1 and at most one rung between any pair of adjacent through strings is the tangle:
\begin{center}
\includegraphics[width=55mm]{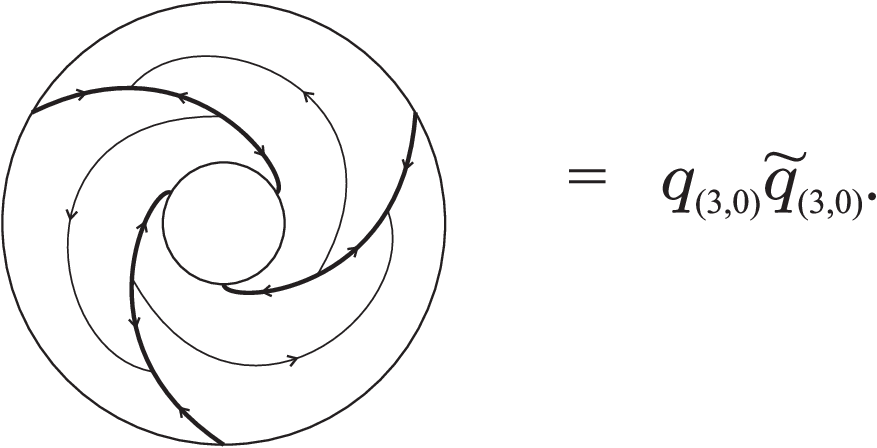}
\end{center}
All tangles with more than one rung between any pair of adjacent through strings are obtained by inserting extra factors of $\widetilde{q}$ or $\widetilde{q}^{\ast}$.

Now suppose that any basis tangle in $A^{(3,0)}$ for which a choice of through strings can be made so that the tangle has rotation number $k$, can be written as $q^i \widetilde{q}^j$, where $i \in \{ 0,1,2 \}$, $j \in \mathbb{N}\cup\{0\}$. The if $T$ can be drawn as a tangle with rotation number $k+1$, $Tq^{\ast}$ can be drawn as a tangle with rotation number $k$, thus $Tq^{\ast} = q^i \widetilde{q}^j$ for some $i \in \mathbb{Z}$, $j \in \mathbb{N}\cup\{0\}$. Then $T = q^{i+1} \widetilde{q}^j$. Hence we have shown that $A^{(3,0)}$ is generated by $q$ and $\widetilde{q}$.
Similarly, $A^{(3k,0)}$ is generated by $q_{(3k,0)}$ and $\widetilde{q}_{(3k,0)}$, for all $k \in \mathbb{N}$.

Then $A^{(3,0)}_0$ is generated by $\varphi_{(3,0)} := q_{(3,0)}\widetilde{q}_{(3,0)}$ and $\widetilde{q}_{(3,0)}^3$. Since $\widetilde{q}_{(3,0)}^3 = \varphi_{(3,0)}^3$, thus $A^{(3,0)}_0$ is generated by $\varphi_{(3,0)}$. It is infinite dimensional, but we can construct a finite-dimensional $A_2$-$ATL^{\overline{0}}_{+^3}$-module $V^{(3,0),\gamma}_{+^3,\overline{0}}$ by letting $\varphi_{(3,0)} (= \varphi_{(3,0)}^{\ast})$ count as some scalar $\gamma \in \mathbb{R}$ in $A^{(3,0)}_0$, i.e. $\varphi_{(3,0)} = \gamma \mathbf{1}_{+^3}$. Then we have an $A_2$-$ATL$-module $V^{(3,0),\gamma} = A_2\mathrm{-}ATL(V^{(3,0),\gamma}_{+^3,\overline{0}})$. Note that for each $a \in \{ 0,1,2 \}$, $V^{(3,0),\gamma}_{+^3,\overline{a}}$ is one-dimensional, $V^{(3,0),\gamma}_{+^3,\overline{a}} = \mathbb{C} q_{(3,0)}^{a}$.

For any two elements $S,T \in V^{(3,0),\gamma}_{\sigma}$, the tangle $T^{\ast}S$ will have three (source) vertices on its outer disc and three (sink) vertices on its inner disc. We use relations K1-K3 on $T^{\ast}S$ to obtain a linear combination $\sum_j c_j (T^{\ast}S)_j$ of tangles $(T^{\ast}S)_l$ which do not contain any closed circles, digons or embedded squares, where $c_j \in \mathbb{C}$. We let $\langle S,T \rangle_l$ be zero if $\mathrm{rank}((T^{\ast}S)_l) < (3,0)$. Otherwise, $(T^{\ast}S)_l$ will be equal to $\varphi_{(3,0)}^{p}$ for some $p=0,1,2,\ldots \;$, and we let $\langle S,T \rangle_l$ be $\gamma^p$. We then define an inner product on $V^{(3,0),\gamma}$ by $\langle S,T \rangle = \sum_j c_j \langle S,T \rangle_j$.

The $A_2$-$ATL$-module $V^{(3,0),\gamma}$ is irreducible by Lemma \ref{Lemma:SU(3)Lemma3.4}.
If the inner product is positive semi-definite we define the Hilbert $A_2$-$ATL$-module $H^{(3,0),\gamma}$ to be the quotient of $V^{(3,0),\gamma}$ by the zero-length vectors with respect to this inner product; otherwise $H^{(3,0),\gamma}$ does not exist.
There is a similar description of modules $H^{(0,3),\gamma}$ of minimum rank $(0,3)$, where there are now three source vertices on the inner disc.

\begin{figure}[htb]
\begin{center}
  \includegraphics[width=120mm]{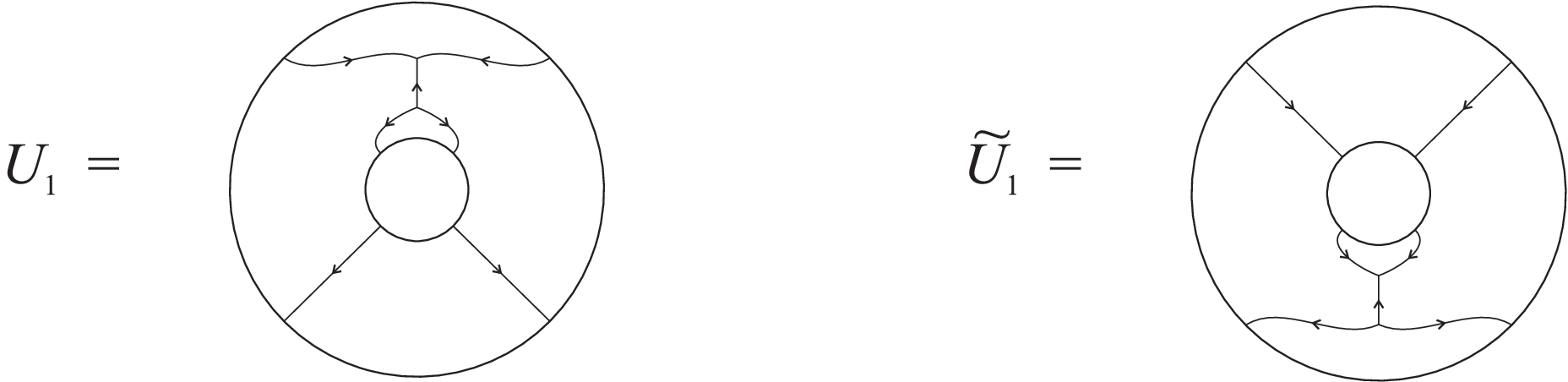}
 \caption{Annular $+^2 -^2$-tangles $U_1$, $\widetilde{U}_1$} \label{fig:UU}
\end{center}
\end{figure}

\begin{figure}[htb]
\begin{center}
  \includegraphics[width=120mm]{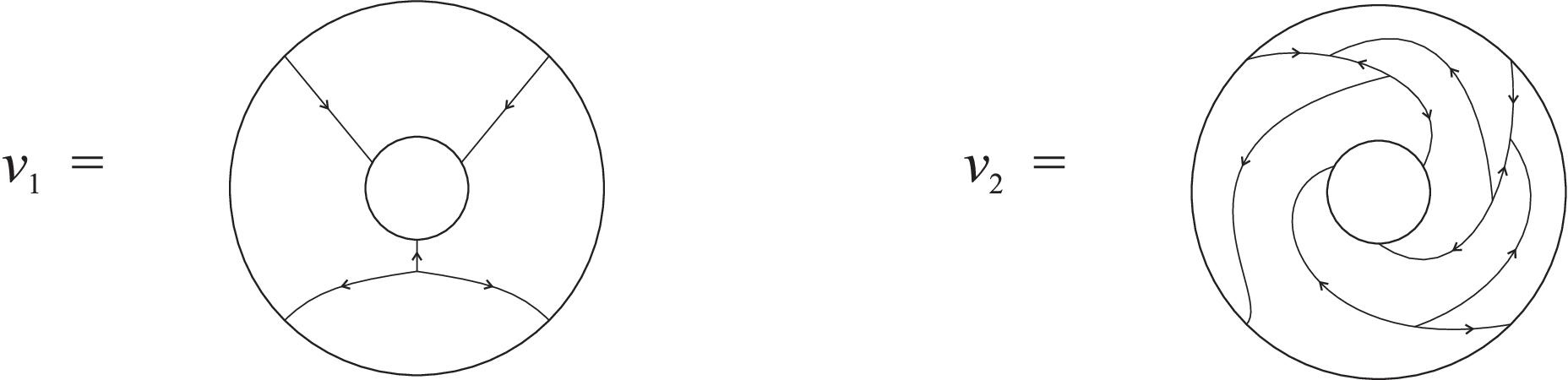}
 \caption{The basis elements $v_1$, $v_2$ of $V^{(3,0),\gamma}_{+^2 -^2}$} \label{fig:v1v2}
\end{center}
\end{figure}

We will now conjecture certain irreducible modules of non-zero weight that the $A_2$-planar algebra $P^{\mathcal{G}}$ for the graphs $\mathcal{E}^{(8)}$ and $\mathcal{D}^{(6)}$ should contain.
Let $U_1, \widetilde{U}_1 \in A_2\textrm{-}ATL_{+^2 -^2}$ be the annular $+^2 -^2$-tangles illustrated in Figure \ref{fig:UU}.
From drawing pictures, it appears that the lowest weight module $V^{(3,0),\gamma}_{+^2 -^2}$ is the span of $v_1, v_2$, illustrated in Figure \ref{fig:v1v2}, where $v_{2} = \varphi_{(2,2)} v_{1}$. These are the only tangles we can find that have rank no smaller than $(3,0)$, do not contain any closed circles, digons or embedded squares, and which cannot be written as a linear combination of tangles of the form $v' \varphi_{(3,0)}^{p}$ for some $p \in \mathbb{N}$, where $v'$ is $v_1$ or $v_2$, and the tangle $\varphi_{(3,0)}^{p}$ is inserted in the inner disc of $v'$.

The action of $A_2$-$ATL_{+^2 -^2}$ on $H^{(3,0),\pm 1}_{+^2 -^2}$ is given explicitly as
$$ \begin{array}{rcllrcll}
\varphi_{(2,2)} v_1 & = & v_2, & & \varphi_{(2,2)} v_2 & = & v_1, & \\
\widetilde{U}_1 v_1 & = & \delta v_1, & & \widetilde{U}_1 v_2 & = & \gamma v_1, & \\
\widetilde{\varphi}_{(2,2)} v_l & = & v_l, & & U_1 v_l & = & 0, & l=1,2.
\end{array}$$
The roles of $U_l$ and $\widetilde{U}_l$ are interchanged for $H^{(0,3),\gamma}$.

We were able to conjecture certain irreducible modules of non-zero weight that the $A_2$-graph planar algebra $P^{\mathcal{G}}$ for the graphs $\mathcal{E}^{(8)}$ and $\mathcal{D}^{(6)}$ should contain, since the action of the rotation $\rho_{(2,2)}$ on the $A_2$-planar algebras for these graphs was much easier to write down than for the other graphs.

For the graph $\mathcal{E}^{(8)}$, its zero-weight irreducible modules are $H^{\beta_{(0,0)}}$, $H^{\beta_{(3,0)}}$, $H^{\beta_{(0,3)}}$ and $H^{\beta_{(2,2)}}$. By computing the inner-products $\langle v_i, v_j \rangle$ of the elements $v_l \in H^{\beta}_{0,1}$ explicitly, and using Mathematica to compute the rank of the matrix $(\langle v_i, v_j \rangle)_{i,j}$, we computed the dimension of $H^{\beta_{(0,0)}}_{+-}$, $H^{\beta_{(3,0)}}_{+-}$, $H^{\beta_{(0,3)}}_{+-}$ and $H^{\beta_{(2,2)}}_{+-}$ and found that $P^{\mathcal{E}^{(8)}}$ did not contain any irreducible modules of lowest weight $(1,1)$.
Similarly, by computing the dimensions of $W = H^{\beta_{(0,0)}}_{+^2 -^2,\overline{0}} \oplus H^{\beta_{(3,0)}}_{+^2 -^2,\overline{0}} \oplus H^{\beta_{(0,3)}}_{+^2 -^2,\overline{0}} \oplus H^{\beta_{(2,2)}}_{+^2 -^2,\overline{0}}$, we find that $\mathrm{dim}(W) = 30$ whilst $\mathrm{dim}(P^{\mathcal{E}^{(8)}}_{+^2 -^2,\overline{0}}) = 36$, so that the dimension of $W^{\bot} \cap P^{\mathcal{E}^{(8)}}_{+^2 -^2,\overline{0}}$ is 6. Then for modules of lowest weight $\leq (2,2)$, we conjecture
$$P^{\mathcal{E}^{(8)}}_{+^2 -^2} = H^{\beta_{(0,0)}}_{+^2 -^2} \oplus H^{\beta_{(3,0)}}_{+^2 -^2} \oplus H^{\beta_{(0,3)}}_{+^2 -^2} \oplus H^{\beta_{(2,2)}}_{+^2 -^2} \oplus H^{(3,0),\varepsilon_1}_{+^2 -^2} \oplus H^{(0,3),\varepsilon_1}_{+^2 -^2} \oplus H^{(2,2),\gamma_1,\varepsilon_2i}_{+^2 -^2} \oplus H^{(2,2),\gamma_2,\varepsilon_3i}_{+^2 -^2},$$
where $\varepsilon_i \in \{ \pm 1 \}$, $i = 1,2,3$, and $\gamma_1, \gamma_2 \in \mathbb{T}$, where the exact values of these six parameters has not yet been determined. This conjecture arises from computing the eigenvalues of the actions of $\rho_{(2,2)}$, $U_1$ and $\widetilde{U}_1$ on $W^{\bot} \cap P^{\mathcal{E}^{(8)}}_{+^2 -^2,\overline{0}}$. Each action is a linear transformation, which we computed by hand, and then computed using Mathematica the eigenvalues of the matrix which gives this linear transformation. These eigenvalues are
\begin{eqnarray}
\rho_{(2,2)}: & & 1 \textrm{ twice, } -1 \textrm{ four times}, \\
U_1, \widetilde{U}_1: & & [4] \alpha \delta^{-2}, \textrm{ once, } 0 \textrm{ five times}.
\end{eqnarray}
The eigenvalues of the actions of these elements on $H^{(2,2),\gamma,\omega}_{+^2 -^2,\overline{a}}$, $H^{(3,0),\gamma}_{+^2 -^2,\overline{a}}$ and $H^{(0,3),\gamma}_{+^2 -^2,\overline{a}}$ are given in the Table \ref{Table:action_on_modules}.

\begin{table}[hbt]
\begin{center}
\begin{tabular}{|c|c|c|c|} \hline
& \multicolumn{3}{|c|}{Eigenvalues of the action of} \\
$A_2$-$ATL$-module & $\rho_{(2,2)}$ & $U_1$ & $\widetilde{U}_1$ \\
\hline\hline $H^{(2,2),\gamma,\omega}_{+^2 -^2,\overline{a}}$ & $\omega^2$ & 0 & 0 \\
\hline $H^{(3,0),\pm 1}_{+^2 -^2,\overline{a}}$ & $1, \;\; -1$ & 0 $\; (\times 2)$ & $[4] \alpha \delta^{-2}, \;\; 0$ \\
\hline $H^{(0,3),\pm 1}_{+^2 -^2,\overline{a}}$ & $1, \;\; -1$ & $[4] \alpha \delta^{-2}, \;\; 0$ & 0 $\; (\times 2)$ \\
\hline $H^{(0,3),\gamma}_{+^2 -^2,\overline{a}}$, $\gamma \neq \pm 1$ & 1 $\, (\times 3), \;\, -1$ $\, (\times 3)$ & 0 $\; (\times 6)$ & $[4] \alpha \delta^{-2} \;\, (\times 3), \;\, 0 \;\, (\times 3)$ \\
\hline $H^{(3,0),\gamma}_{+^2 -^2,\overline{a}}$, $\gamma \neq \pm 1$ & 1 $\, (\times 3), \;\, -1$ $\, (\times 3)$ & $[4] \alpha \delta^{-2} \;\, (\times 3), \;\, 0 \;\, (\times 3)$ & 0 $\; (\times 6)$ \\
\hline
\end{tabular}\\
\caption{The eigenvalues of the actions of $\rho_{(2,2)}$, $U_1$, $\widetilde{U}_1$ on $H^{(2,2),\gamma,\omega}_{+^2 -^2,\overline{a}}$, $H^{(3,0),\gamma}_{+^2 -^2,\overline{a}}$, $H^{(0,3),\gamma}_{+^2 -^2,\overline{a}}$.} \label{Table:action_on_modules}
\end{center}
\end{table}

Then we see that $W^{\bot} \cap P^{\mathcal{E}^{(8)}}_{+^2 -^2}$ should contain one copy of both of $H^{(3,0),\varepsilon_1}_{+^2 -^2}$ and $H^{(0,3),\varepsilon_1'}_{+^2 -^2}$, $\varepsilon_1, \varepsilon_1' \in \{ \pm 1 \}$, and since $P^{\mathcal{E}^{(8)}}$ is invariant under conjugation of the graph $\mathcal{E}^{(8)}$, we should have $\varepsilon_1 = \varepsilon_1'$. Then we need two rank $(2,2)$ modules $H^{(2,2),\gamma_1,\omega}_{+^2 -^2}$, $H^{(2,2),\gamma_2,\omega}_{+^2 -^2}$ such that the action of $\rho_{(2,2)}$ on both has an eigenvalue $\omega^2 = -1$, i.e. $\omega = \pm i$. Since $P^{\mathcal{E}^{(8)}}$ is invariant under complex conjugation, we would either have $\gamma_1, \gamma_2 \in \mathbb{R}$ or else $\gamma_1 = \overline{\gamma}_2$. However, to determine the exact values of $\varepsilon_i$, $i=1,2,3$, and $\gamma_1$, $\gamma_2$, we would need to consider the action of $\varphi_{(2,2)}$ on $W^{\bot} \cap P^{\mathcal{E}^{(8)}}_{+^2 -^2}$, the computation of which is extremely tedious. So we have
$$P^{\mathcal{E}^{(8)}} \supset H^{\beta_{(0,0)}} \oplus H^{\beta_{(3,0)}} \oplus H^{\beta_{(0,3)}} \oplus H^{\beta_{(2,2)}} \oplus H^{(3,0),\varepsilon_1} \oplus H^{(0,3),\varepsilon_1} \oplus H^{(2,2),\gamma_1,\varepsilon_2i} \oplus H^{(2,2),\gamma_2,\varepsilon_3i}.$$

Similarly for the graph $\mathcal{D}^{(6)}$, $P^{\mathcal{D}^{(6)}}$ contains no irreducible modules of lowest weight $(1,1)$. Computing the dimensions of $P^{\mathcal{D}^{(6)}}_{+^2 -^2}$ and $W = H^{\beta_{(0,0)}}_{+^2 -^2,\overline{0}} \oplus H^{0,\overline{0}}_{+^2 -^2,\overline{0}}$ as for the $\mathcal{E}^{(8)}$ case, we find $\mathrm{dim}(P^{\mathcal{D}^{(6)}}_{+^2 -^2,\overline{0}}) = 16$ and $\mathrm{dim}(W) = 14$. Then the dimension of $W^{\bot} \cap P^{\mathcal{D}^{(6)}}_{+^2 -^2}$ is 2, and hence $P^{\mathcal{D}^{(6)}}_{+^2 -^2}$ must either contain one copy of $H^{(3,0),\gamma}_{+^2 -^2}$ or else $H^{(2,2),\gamma_1,\omega_1}_{+^2 -^2} \oplus H^{(2,2),\gamma_2,\omega_2}_{+^2 -^2}$. By considering the action of $\rho_{(2,2)}$ on $W^{\bot} \cap P^{\mathcal{D}^{(6)}}_{+^2 -^2}$, we have the eigenvalue 1 twice. Then $W = H^{(2,2),\gamma_1,\omega_1}_{+^2 -^2} \oplus H^{(2,2),\gamma_2,\omega_2}_{+^2 -^2}$, where $\omega_i^2 = 1$, $i=1,2$. Then we see that
$$P^{\mathcal{D}^{(6)}} \supset H^{\beta_{(0,0)}} \oplus H^{0,\overline{0}} \oplus H^{(2,2),\gamma_1,\varepsilon_1} \oplus H^{(2,2),\gamma_2,\varepsilon_2},$$
where $\varepsilon_1, \varepsilon_2 \in \{ \pm 1 \}$, and either $\gamma_1, \gamma_2 \in \mathbb{R}$ or else $\gamma_1 = \overline{\gamma}_2$. Again, to determine the values of $\varepsilon_i, \gamma_i$, $i=1,2$, explicitly requires considering the eigenvalues of the action of $\varphi_{(2,2)}$ on $W^{\bot} \cap P^{\mathcal{D}^{(6)}}_{+^2 -^2}$.

\paragraph{Acknowledgements}

This paper is based on work in \cite{pugh:2008}. The first author was partially supported by the EU-NCG network in Non-Commutative Geometry MRTN-CT-2006-031962, and the second author was supported by a scholarship from the School of Mathematics, Cardiff University.

\end{document}